\documentclass[review,onefignum,onetabnum]{siamart171218}

\usepackage{lipsum}
\usepackage{amsfonts}
\usepackage{graphicx}
\usepackage{epstopdf}
\usepackage{algorithmic}
\usepackage{amssymb}
\usepackage[english]{babel}
\usepackage{array}
\usepackage{adjustbox}
\usepackage{relsize}
\usepackage{spverbatim}
\usepackage{listings}
\usepackage{mathrsfs}
\usepackage[skip=4pt plus1pt, indent=20pt]{parskip}

\newcolumntype{L}{>{$}l<{$}}

\ifpdf
  \DeclareGraphicsExtensions{.eps,.pdf,.png,.jpg}
\else
  \DeclareGraphicsExtensions{.eps}
\fi

\newsiamremark{hypothesis}{Hypothesis}
\crefname{hypothesis}{Hypothesis}{Hypotheses}
\newsiamthm{claim}{Claim}
\newtheorem{remark}{Remark}
\nolinenumbers

\headers{Zeros of Bessel functions}{T. M. Dunster}

\title{Error bounds for a uniform asymptotic approximation of the zeros of the Bessel function $J_{\nu}(\MakeLowercase{x})$}

\author{T. M. Dunster\thanks{Department of Mathematics and Statistics, San Diego State University, 5500 Campanile Drive, San Diego, CA 92182-7720, USA. 
  (\email{mdunster@sdsu.edu}, \url{https://tmdunster.sdsu.edu}).}
  }

\usepackage{amsopn}

\makeatletter
\newcommand*{\addFileDependency}[1]{
  \typeout{(#1)}
  \@addtofilelist{#1}
  \IfFileExists{#1}{}{\typeout{No file #1.}}
}
\makeatother

\nolinenumbers

\begin{document}

\maketitle

\begin{abstract}
A recent asymptotic expansion for the positive zeros $x=j_{\nu,m}$ ($m=1,2,3,\ldots$) of the Bessel function of the first kind $J_{\nu}(x)$ is studied, where the order $\nu$ is positive. Unlike previous well-known expansions in the literature, this is uniformly valid for one or both $m$ and $\nu$ unbounded, namely $m=1,2,3,\ldots$ and $1 \leq \nu < \infty$. Explicit and simple lower and upper error bounds are derived for the difference between $j_{\nu,m}$ and the first three terms of the expansion. The bounds are sharp in the sense they are close to the value of the fourth term of the expansion (i.e. the first neglected term).
\end{abstract}



\section{Introduction and main result}

We consider the positive zeros $x=j_{\nu,m}$ ($m=1,2,3,\ldots$) of the Bessel function of the first kind $J_{\nu}(x)$, ordered by increasing values. Recently in \cite{Dunster:2024:AZB}, and based on the classic 1954 paper of Olver \cite{Olver:1954:AEB}, asymptotic expansions were constructed by the present author of the form
\begin{equation}
\label{eq1.2}
j_{\nu,m} \sim 
\nu \sum_{s=0}^{\infty}\frac{z_{m,s}}{\nu^{2s}},
\end{equation}
uniformly for $m=1,2,3,\ldots$ and $\nu > 0$. The coefficients $z_{m,s}$ ($s=1,2,3,\ldots$) were shown to be rational functions of the leading term $z_{m,0}$ and two other readily computable variables (described below). Apart from $z_{m,0}$ these coefficients vanish as $m \to \infty$, thus providing a powerful and simple uniform asymptotic expansion for one or both of $m$ or $\nu$ large. The first four coefficients are given below, and the rest can be evaluated via a recurrence relation given in \cite{Dunster:2024:AZB}.

The purpose of this paper is to provide error bounds for a truncated version of (\ref{eq1.2}). Sharp, simple and useful error bounds for the well-known expansions \cite[Eqs. 10.21.19 and 10.21.32]{NIST:DLMF}, where only one of $m$ and $\nu$ is permitted to be large, were proven by Nemes \cite{Nemes:2021:PTC} ($m \in \mathbb{Z}^{+}$ and $\nu \in (-\tfrac{1}{2},\tfrac{1}{2}) $), and Qu and Wong \cite{Qu:1999:BPU} ($\nu \in (0, \infty)$ with $m$ bounded). See also references therein for earlier results in both cases. For discussions on the significance of explicit error bounds in asymptotic approximations see \cite{Olver:1980:AEB} and \cite{Wong:1980:EBA}.

Our main result is given by \cref{thm:main} below, which we prove in \cref{sec2}. The proof primarily utilises \cref{thm:AiryError} which is proven in \cref{sec3}, \cref{thm:e_bold} which is proven in \cref{sec4}, as well as \cref{thm:hethcote,thm:QuWong}, the proofs of which are given in the papers by Hethcote \cite{Hethcote:1970:EBA}, and Qu and Wong \cite{Qu:1999:BPU}, respectively. In this and the next sections we also state a number of lemmas which we shall use, and the proofs of these are deferred to \cref{secA}.

Returning to the expansion (\ref{eq1.2}), in order to describe the coefficients we first define
\begin{equation}
\label{eq1.3}
\frac{2}{3}(-\zeta)^{3/2}
=\int_{1}^{z}\frac{\left(t^2-1\right)^{1/2}}{t}dt
=\left(z^{2}-1\right)^{1/2}
-\mathrm{arcsec}(z)
\quad  (1 \leq z < \infty).
\end{equation}
This is the Liouville variable that is used in the standard Airy function approximations of Bessel functions (see \cite[Chap. 11]{Olver:1997:ASF}). Principal branches are taken in (\ref{eq1.3}) so that the interval $1 \leq z < \infty$ is mapped to $-\infty < \zeta \leq 0$. Further, $\zeta$ is analytic at $z=1$, which is a turning point of the equation which is satisfied by $J_{\nu}(\nu z)$, namely
\begin{equation}
\label{eq1.3a}
\frac{d^{2} w}{dz^{2}}=
\left\{\nu^{2}\frac{1-z^{2}}{z^{2}} 
-\frac{1}{4z^{2}}\right\} w.
\end{equation}

From (\ref{eq1.3}) it is straightforward to show that $\zeta \to 0$ as $z \to 1$ such that
\begin{equation}
\label{eq1.3b}
\zeta = -2^{1/3}(z-1)
+\tfrac{3}{10}2^{1/3}(z-1)^2
-\tfrac{32}{175}2^{1/3}(z-1)^3
+\mathcal{O}\left\{(z-1)^4\right\},
\end{equation}
and $\zeta \to -\infty$ as $z \to \infty$ such that
\begin{equation}
\label{eq1.3c}
\zeta = -\tfrac{1}{4} 12^{2/3}z^{2/3}
\left\{1-\tfrac{1}{3}z^{-1}
+\mathcal{O}\left(z^{-2}\right)\right\}.
\end{equation}
Both (\ref{eq1.3b}) and (\ref{eq1.3c}) are repeatedly differentiable.

Using (\ref{eq1.3}) one can readily verify that $\zeta$ is monotonically decreasing for $1 < z < \infty$. We also need similar results for its first few derivatives, as given in the following.
\begin{lemma}
\label{lem:zeta}
For $1 \leq z < \infty$, $-\zeta'$, $\zeta''$ and $-\zeta'''$ are positive and decrease monotonically to zero.
\end{lemma}

Next, for $m=1,2,3,\ldots$ and $\nu > 0$, the leading coefficient $z_{m,0}$ of the expansion (\ref{eq1.2}) is the unique value lying in $(1,\infty)$ that satisfies
\begin{equation}
\label{eq1.5}
\zeta_{m,0}:=\zeta(z_{m,0})
=\nu^{-2/3}\mathrm{a}_{m},
\end{equation}
where $0>\mathrm{a}_{1}>\mathrm{a}_{2}>\mathrm{a}_{3}>\ldots$ are the negative zeros of the Airy function $\mathrm{Ai}(x)$ (see \cite[Sect. 9.9]{NIST:DLMF}). Thus from (\ref{eq1.3}) $z_{m,0}$ is the solution of the implicit equation
\begin{equation}
\label{eq1.4}
\begin{split}
\int_{1}^{z_{m,0}}
\frac{\left(t^{2}-1\right)^{1/2}}{t}\,dt
& =\left(z_{m,0}^{2}-1\right)^{1/2}
-\mathrm{arcsec}\left( z_{m,0} \right)
\\ &
= \frac{2}{3}\left|\zeta_{m,0}\right|^{3/2}
=\frac{2}{3 \nu}\left|\mathrm{a}_{m}\right|^{3/2}.
\end{split}
\end{equation}

From (\ref{eq1.5}) and \cite[\S 9.9(iv)]{NIST:DLMF}
\begin{equation}
\label{eq1.5a}
\zeta_{m,0}=\nu^{-2/3}
\mathrm{a}_{m,0}\left\{
1+\mathcal{O}\left(m^{-2}\right)\right\}
\quad (m \to \infty),
\end{equation}
where
\begin{equation}
\label{eq1.50}
\mathrm{a}_{m,0}=-\left\{
\tfrac{3}{8}\pi(4m-1)\right\}^{2/3}.
\end{equation}
For $1 \leq \nu < \infty$ and $m=1,2,3,\ldots$ we have $-\infty < \zeta_{m,0} <0$, and $1<z_{m,0}<\infty$ with
\begin{equation}
\label{eq1.6}
z_{m,0} = 1 -2^{-1/3}\zeta_{m,0}
+\mathcal{O}\left(\zeta_{m,0}^{2}\right)
\quad (\zeta_{m,0} \to 0^{-}),
\end{equation}
and
\begin{equation}
\label{eq1.7}
z_{m,0} = \tfrac{2}{3}\left|\zeta_{m,0}\right|^{3/2}
+\tfrac{1}{2}\pi+\mathcal{O}\left(
\left|\zeta_{m,0}\right|^{-3/2}
\right)
\quad (\zeta_{m,0} \to -\infty).
\end{equation}
Note from (\ref{eq1.4}) that $z_{m,0} \to 1$ ($\zeta_{m,0} \to 0$) occurs when $\nu^{-2/3}\mathrm{a}_{m} \to 0$, such as $\nu \to \infty$ with $m$ being fixed, and $z_{m,0} \to \infty$ ($\zeta_{m,0} \to -\infty$) occurs when $\nu^{-2/3}|\mathrm{a}_{m}| \to \infty$, such as $m \to \infty$ with $\nu$ fixed.

As shown in \cite{Dunster:2024:AZB}, the later coefficients in (\ref{eq1.2}) are rational functions of the first coefficient $z_{m,0}$, as well as $\zeta_{m,0}$ and $\sigma_{m,0}:=\sigma(z_{m,0})$, where
\begin{equation}
\label{eq1.8}
\sigma(z)
=\left(\frac{\zeta}{ 1-z^{2}}\right)^{1/2}.
\end{equation}
This function is analytic at $z=1$, and in particular from (\ref{eq1.3b})
\begin{equation}
\label{eq1.9}
\sigma(z)
=\frac{1}{2^{1/3}}
-\frac{{2}^{2/3}}{5}(z-1)
+\mathcal{O}\left\{(z-1)^2\right\}
\quad (z \to 1),
\end{equation}
and also, from (\ref{eq1.3c}),
\begin{equation}
\label{eq1.9b}
\sigma(z)
=\frac{12^{1/3}}{2z^{2/3}}
\left\{1+\mathcal{O}\left(
\frac{1}{z}\right)\right\}
\quad (z \to \infty).
\end{equation}

From (\ref{eq1.3}) and (\ref{eq1.8}) the following derivatives are obtained and will later be used
\begin{equation}
\label{eq1.35a}
\frac{d \zeta}{d z}
= -\frac{1}{z \sigma},
\quad
\frac{d \sigma}{dz}
=\frac {2 z^{2} \sigma^{3}-1}{2 z\zeta}.
\end{equation}

Now the coefficients in (\ref{eq1.2}) are given by $z_{m,j}=\hat{z}_{j}(z_{m,0})$, where $\hat{z}_{j}=\hat{z}_{j}(z)$ for $j=1,2$ are given explicitly by
\begin{equation}
\label{eq1.10}
\hat{z}_{1}=
\frac{z \sigma}{48  
\zeta^{2}}
\left\{ 10\sigma^3
- 6\sigma \zeta - 5 \right\},
\end{equation}
and
\begin{equation}
\label{eq1.11}
\begin{split}
\hat{z}_{2}  = &
\frac{z \sigma}{46080
\zeta^{5}}\left\{
200\sigma^{9}\left( 35z^{2}
+221 \right) 
-80\sigma^{7}  \zeta\left( 75z^{2}+982\right) \right.
\\
& -4000\sigma^{6}  z^{2}
+24\sigma^{5}  \zeta^{2}\left( 45 z^{2}+1543\right)
+200\sigma^{4}  \zeta
\left(6 z^{2}-5 \right) 
 \\
& \left. +600\sigma^{2}  \zeta^{2}
+10 \sigma^{3}\left(25 z^{2} 
-264 \zeta^{3}\right) 
+250  \sigma \zeta-5525
\right\}.
\end{split}
\end{equation}

\begin{lemma}
\label{lem:zmjMonotonic}
For $1 \leq z < \infty$  $\hat{z}_{1}(z)$, $-\hat{z}_{2}(z)$ and $\hat{z}_{3}(z)$ are positive and monotonically decreasing to 0. Moreover
\begin{equation}
\label{eq1.20}
\hat{z}_{1}(z)+\hat{z}_{2}(z)
>\hat{z}_{1}(z)+2 \hat{z}_{2}(z)>0.
\end{equation}
\end{lemma}

Our focus is the error in the truncated approximation
\begin{equation}
\label{eq1.12}
\frac{j_{\nu,m}}{\nu} \approx 
\mathbf{z}_{\nu,m} := z_{m,0}+
\frac{z_{m,1}}{\nu^{2}}
+\frac{z_{m,2}}{\nu^{4}}.
\end{equation}

\begin{lemma}
\label{lem:bold(z[m])Bound}
For $m=1,2,3,\ldots$, $1 \leq \nu < \infty$, and hence $1 < z_{m,0} < \infty$, we have
\begin{equation}
\label{eq1.21}
z_{m,0}< \mathbf{z}_{\nu,m} <
z_{m,0}+\tfrac{1}{76}.
\end{equation}
\end{lemma}

In our error bound for (\ref{eq1.12}) the next coefficient $z_{m,3}$ from the full expansion (\ref{eq1.2}) appears, which is given by $z_{m,3}=\hat{z}_{3}(z_{m,0})$ where $\hat{z}_{3}=\hat{z}_{3}(z)$ has the explicit representation
\begin{equation}
\begin{split}
\label{eq1.13}
\hat{z}_{3}=&
\frac{z\sigma}{92897280 \zeta^{8}}
\left\{
28000 \sigma^{15}\left( 650z^{4}+8619z^{2}+89451 \right) \right.
\\ &
-5600\sigma^{13}\zeta \left( 3960z^{4}+86491z^{2}+1338450 \right) 
-420000z^{2} \sigma^{12}\left( 35z^{2}+221 \right) 
\\ &
+1120  \sigma^{11} \zeta^{2}
\left( 7290z^{4}+282717z^{2}+7233983 \right)
\\ &
+8400 \sigma^{10} \zeta
\left( 1200z^{4}+15187z^{2}-1105 \right) 
\\ &
-672\sigma^{9}\left(1260z^{4}\zeta^{3}
+107859 z^{2}\zeta^{3}-4375z^{4}+
5551275\zeta^{3} \right) 
\\ &
-3360 \sigma^{8} \zeta^{2} 
\left( 405z^{4}+12762z^{2}-4910 \right) 
\\ &
+24 \sigma^{7}\zeta \left(136080 
z^{2}\zeta^{3}-26250 z^{4}
+26578824 \zeta^{3}+74375 z^{2} \right) 
\\ &
+560 \sigma^{6}\left( 2745   z^{2}\zeta^{3}-125 z^{4}
-13887\zeta^{3}-16575z^{2} \right)
\\ &
-168  \sigma^{5} \zeta^{2} 
\left( 111564\zeta^{3}+3375 z^{2}-
625 \right)
\\ &
+700 \sigma^{4}\zeta \left( 792 \sigma^{3} \zeta^{3}+3853z^{2}-3315
\right) -5250\sigma^{3}
\left( 12\zeta^{3}-221z^{2} \right) 
\\ &
\left. 
+1374800\sigma^{2}\zeta^{2}+1160250
\sigma\zeta-78269625
\right\}.
\end{split}
\end{equation}

As $z \to 1$ ($\zeta \to 0$) we have from (\ref{eq1.3b}), (\ref{eq1.9}), (\ref{eq1.10}), (\ref{eq1.11}) and (\ref{eq1.13})
\begin{equation}
\label{eq1.14}
\hat{z}_{1}(z)=
\tfrac{1}{70} 
-\tfrac{1}{3150} (z-1)
+\mathcal{O}\left\{(z-1)^{2}\right\}),
\end{equation}
\begin{equation}
\label{eq1.15}
\hat{z}_{2}(z)=
-\tfrac{3781}{3185000} 
+\tfrac{5750429}{9029475000} (z-1)
+\mathcal{O}\left\{(z-1)^{2}\right\}),
\end{equation}
and
\begin{equation}
\label{eq1.16}
\hat{z}_{3}(z) =
\tfrac{722735647}{1630879250000} 
-\tfrac{1574852287133}{3566732919750000} 
(z-1)
+\mathcal{O}\left\{(z-1)^{2}\right\}),
\end{equation}
and from (\ref{eq1.3c}), (\ref{eq1.9b}), (\ref{eq1.10}), (\ref{eq1.11}) and (\ref{eq1.13}) as $z \to \infty$ ($\zeta \to -\infty$)
\begin{equation}
\label{eq1.17}
\hat{z}_{1}(z) =
\tfrac{1}{18} z^{-1}
+\mathcal{O}\left(z^{-2}\right),
\end{equation}
\begin{equation}
\label{eq1.18}
\hat{z}_{2}(z) =
-\tfrac{71}{1944 } z^{-3}
+\mathcal{O}\left(z^{-4}\right),
\end{equation}
and
\begin{equation}
\label{eq1.19}
\hat{z}_{3}(z) =
\tfrac{6673}{58320} z^{-5}
+\mathcal{O}\left(z^{-6}\right),
\end{equation}
with all six of the above being repeatedly differentiable.

\begin{remark}
The proof in \cref{secA} of \cref{lem:zmjMonotonic,lem:bold(z[m])Bound} and a number of the other results relies on numerical evaluation of certain explicitly-given functions of $z$. For example, on examining (\ref{eq1.13}) we see that proving positivity and monotonicity of $\hat{z}_{3}(z)$ by analytical means would be extremely challenging, and even if possible almost certainly very long and tedious. On the other hand, it and the other functions studied are explicitly given, and are polynomials of $z$ and the elementary functions $1/\zeta$ and $\sigma$. Thus, for example, to prove $\hat{z}_{3}(z)$ is positive and decreasing for $1\leq z<\infty$ we instead first evaluate $-d \hat{z}_{3}(z)/dz$ explicitly via (\ref{eq1.35a}), followed by a change of variable $z=1/(1-v)$ in this expression. Then a plot of a suitably scaled version of this function verifies it is positive for $0 \leq v <1$ (with $v=1$ being a removable singularity), thus proving the monotonicity of $\hat{z}_{3}(z)$. Positivity is  confirmed by integration of its derivative along with its behaviour (\ref{eq1.19}) at infinity.
\end{remark}

The main result of this paper reads as follows.
\begin{theorem}
\label{thm:main}
For $m=1,2,3,\ldots$ and $1 \leq \nu < \infty$
\begin{equation}
\label{eq1.72}
0  <\frac{z_{m,3}}{\nu^{6}} \left\{
0.969746 
- \frac{\chi_{m}}{\nu^{5/3}}
\right\}
<\frac{j_{\nu,m}}{\nu} - \mathbf{z}_{\nu,m}
<\frac{z_{m,3}}{\nu^{6}} \left\{
1.013023
+ \frac{\chi_{m}}{\nu^{5/3}}\right\},
\end{equation}
where $\mathbf{z}_{\nu,m}$ is defined by (\ref{eq1.12}), and
\begin{equation}
\label{eq1.71a}
\chi_{m}= \frac{2.297225\,
\sigma(1.01354\,z_{m,0})}
{\left|\mathrm{a}_{m,0}\right|^{1/2}},
\end{equation}
in which $\mathrm{a}_{m,0}$ is given by (\ref{eq1.50}), and $\sigma(z)$ by (\ref{eq1.3}) and (\ref{eq1.8}).
\end{theorem}

\begin{remark}
From (\ref{eq1.5a}), (\ref{eq1.50}), (\ref{eq1.7}) and (\ref{eq1.9b}) it is seen that $\nu^{-5/3}\chi_{m}=\mathcal{O}(\min\{m^{-1} \nu^{-1},m^{-1/3}\nu^{-5/3}\})$ uniformly for $m=1,2,3,\ldots$ and $1 \leq \nu < \infty$. Note also from (\ref{eq1.19}) that $z_{m,3}=\mathcal{O}(\{z_{m,0}\}^{-5})$ as $z_{m,0} \to \infty$ (which as we indicated above occurs, for example, when $m \to \infty$ with $\nu$ fixed).
\end{remark}

\begin{remark}
The constants in these bounds were derived using certain estimates that involve our smallest assumed value $\nu=1$. They could be sharpened, that is with values of the two constants in (\ref{eq1.72}) closer to the value $1$, as well as a smaller leading constant in (\ref{eq1.71a}), if our smallest value of $\nu$ was taken to be larger; the structure of the bounds would remain the same. We did not attempt to include further $\nu$ dependence to sharpen the bounds more generally for larger $\nu$ as this would be prohibitively complicated.
\end{remark}

In \cite{Dunster:2024:AZB} the expansion (\ref{eq1.2}) was established via approximations of Bessel functions involving Airy functions having an argument $\nu^{2/3}\mathcal{Z}(\nu,z)$. Here $\mathcal{Z}(\nu,z)$ is a certain function that is analytic at the turning point $z=1$ and was shown to possess an asymptotic expansion of the form 
\begin{equation}
\label{eq1.21a}
\mathcal{Z}(\nu,z) \sim \zeta
+\sum_{s=1}^{\infty}
\frac{\Upsilon_{s}(z)}{\nu^{2s}}
\quad (\nu \to \infty),
\end{equation}
where $\Upsilon_{s}(z)$ ($s=1,2,3,\ldots$) are recursively given polynomials in $z$, $1/\zeta$ and $\sigma$, and also have a removable singularity at $z=1$ ($\zeta=0$). In deriving error bounds in the present paper we shall use a truncated form of this expansion, namely
\begin{equation}
\label{eq1.22}
\mathcal{Z}_{3}(\nu,z)
= \zeta
+\eta(\nu,z),
\end{equation}
where
\begin{equation}
\label{eq1.23}
\eta(\nu,z)
= \frac{\Upsilon_{1}(z)}{\nu^{2}}
+\frac{\Upsilon_{2}(z)}{\nu^{4}}
+\frac{\Upsilon_{3}(z)}{\nu^{6}}.
\end{equation}
In this, the coefficients are explicitly given by
\begin{equation}
\label{eq1.24}
\Upsilon_{1}(z)
=\frac{1}{ 48\,\zeta^{2}}
\left\{
10\,\sigma^{3}-6\,\sigma\,\zeta-5
\right\},
\end{equation}
\begin{equation}
\begin{split}
\label{eq1.25}
\Upsilon_{2}(z)
= & \frac{1}{ 11520\,\zeta^{5}}
\left\{
11050\,\sigma^{9}-19890\,\sigma^{7}\zeta
+9558\,\sigma^{5}\zeta^{2}
-125\,\sigma^{6}+150\,\sigma^{4}\zeta
\right.
\\ &
\left.
-45\,\sigma^{2}\zeta^{2}-250 \left( 3\,\zeta^{3}-2 \right) {
\sigma}^{3}-300\,\sigma\,\zeta-1600
\right\},
\end{split}
\end{equation}
and
\begin{equation}
\label{eq1.26}
\begin{split}
\Upsilon_{3}(z)
= & \frac{1}{5806080 \,\zeta^{8}}
\left\{
156539250\,\sigma^{15}-469617750\,\sigma^{13}\zeta+509154660\,{
\sigma}^{11}\zeta^{2}
\right.
\\ &
-580125\,\sigma^{12}
+1392300\,\sigma^{10}\zeta-1128330
\,\sigma^{8}\zeta^{2}
-140\left( 1681389\,\zeta^{3}-8350 \right) 
\sigma^{9}
\\ &
+90\left( 450441\,\zeta^{3}-23380 \right) \sigma^{7}\zeta-378\left( 3219\,\zeta^{3}
-2680 \right) \sigma^{5}\zeta^{2}
\\ &
+147 \left( 2316\,\zeta^{3}-625 \right) \sigma^{6}
-7875\left( 3\,\zeta^{3}-14 \right)
\sigma^{4}\zeta
-33075\sigma^{2}\zeta^{2}
\\ &
\left.
-6720 \left( 12\,\zeta^{3}-125 \right) 
\sigma^{3}
-504000\,\sigma\,\zeta-5398750
\right\}.
\end{split}
\end{equation}

As $z \to 1$ from (\ref{eq1.3b}) and (\ref{eq1.9})
\begin{equation}
\label{eq1.27}
2^{-1/3}\Upsilon_{1}(z)=
\tfrac{1}{70} 
-\tfrac{2}{225}(z-1)
+\tfrac{953}{202125}(z-1)^2
-\tfrac{17942}{7882875}(z-1)^3
+\mathcal{O}\left\{(z-1)^{4}\right\},
\end{equation}
\begin{equation}
\label{eq1.28}
\begin{split}
2^{-1/3} \Upsilon_{2}(z)=&
-\tfrac{82}{73125} 
+\tfrac{253328}{191008125} (z-1)
-\tfrac{60940232}{70354659375} (z-1)^2
\\ &
+\tfrac{48029412512}{171548111109375} (z-1)^3
+\mathcal{O}\left\{(z-1)^{4}\right\},
\end{split}
\end{equation}
and
\begin{equation}
\label{eq1.28a}
\begin{split}
2^{-1/3} \Upsilon_{3}(z)= &
\tfrac{53780996}{127020403125} 
-\tfrac{19655910112}{28493637046875} (z-1)
+\tfrac{75866644041303056}{145415329603528359375} (z-1)^2
\\ &
-\tfrac{17817798902996416}{198293631277538671875} (z-1)^3
+\mathcal{O}\left\{(z-1)^{4}\right\},
\end{split}
\end{equation}
and as $z \to \infty$, from (\ref{eq1.3c}) and (\ref{eq1.9b}),
\begin{equation}
\label{eq1.29}
\Upsilon_{1}(z)=
\tfrac{1}{108}12^{2/3}z^{-4/3}
\left\{1+\mathcal{O}\left(z^{-1}\right)\right\},
\end{equation}
\begin{equation}
\label{eq1.30}
\Upsilon_{2}(z)=
-\tfrac{4}{729}12^{2/3}z^{-10/3}
\left\{1+\mathcal{O}\left(z^{-1}\right)\right\},
\end{equation}
and
\begin{equation}
\label{eq1.31}
\Upsilon_{3}(z)=
\tfrac{5218}{295245}12^{2/3}z^{-16/3}
\left\{1+\mathcal{O}\left(z^{-1}\right)\right\}.
\end{equation}
All six of these expansions are repeatedly differentiable. 

The following will be used in the subsequent sections.
\begin{lemma}
\label{lem:UpsilonMonotonic}
$\Upsilon_{1}(z)$, $-\Upsilon_{2}(z)$, $\Upsilon_{3}(z)$, $-\Upsilon_{1}'(z)$, $\Upsilon_{2}'(z)$, $-\Upsilon_{3}'(z)$,  $\Upsilon_{1}''(z)$, $\eta(1,z)$, $-\eta'(1,z)$, $-z\eta'(1,z)$ and $z \sigma(z) \eta(1,z)$ are positive for $1 \leq z < \infty$, and are strictly decreasing asymptotically to zero. Moreover
\begin{equation}
\label{eqA.18}
\Upsilon_{2}^{2}(z)
-\tfrac{64}{25}\Upsilon_{1}(z)\Upsilon_{3}(z)<0,
\end{equation}
and
\begin{equation}
\label{eqA.22}
\Upsilon_{2}'^{2}(z)
-3\Upsilon_{1}'(z)\Upsilon_{3}'(z)<0.
\end{equation}
\end{lemma}

\begin{lemma}
\label{lem:SumUpsilon}
For $1 \leq z < \infty$ and $1 \leq \nu < \infty$
\begin{equation}
\label{eq1.32}
0 <  \eta(\nu,z)
\leq \eta(1,z) \leq \eta(1,1)
= \tfrac {44873962351}{3302530481250} 2^{1/3},
\end{equation}
\begin{equation}
\label{eq1.33}
0 < \nu^{2/3} \eta(\nu,z)
\leq  \eta(1,1),
\end{equation}
and
\begin{equation}
\label{eq1.34}
-\tfrac {21376589042}{2590330640625} 2^{1/3}
= \eta'(1,1) \leq  \eta'(1,z)
\leq  \eta'(\nu,z) < 0.
\end{equation}
\end{lemma}

\begin{lemma}
\label{lem:ScriptZ}
For $1 \leq z < \infty$ and $1 \leq \nu < \infty$
\begin{equation}
\label{eq1.35}
\mathcal{Z}_{3}'(\nu,z)
:= \frac{\partial \mathcal{Z}_{3}(\nu,z) }
{\partial z}
< -\frac{1}{z \sigma(z)}<0,
\end{equation}
and $\mathcal{Z}_{3}(\nu,z)<0$ for $z_{m,0} \leq z < \infty$ ($-\infty<\zeta \leq \zeta_{m,0}$) and $1 \leq \nu < \infty$.
\end{lemma}

\emph{Unless otherwise stated, here and throughout this paper primes and dots represent derivatives with respect to $z$ and $\zeta$, respectively}. Thus for example from (\ref{eq1.35a}) and the chain rule $\dot{\eta}=-z \sigma \eta'$. With this notation we have the following, the proof of which is given in \cref{sec3}.

\begin{theorem}
\label{thm:AiryError}
Let
\begin{equation}
\label{eq1.38}
\mathcal{A}i(\nu,z)
=\left(1+\dot{\eta}\right)^{-1/2}
\mathrm{Ai}\left( \nu^{2/3}\mathcal{Z}_{3}(\nu,z)\right),
\end{equation}
and
\begin{equation}
\label{eq1.39}
M(x)=\left\{\mathrm{Ai}^{2}(x)
+\mathrm{Bi}^{2}(x)\right\}^{1/2}
\quad (-\infty < x < \infty).
\end{equation}
Then 
\begin{equation}
\label{eq1.40}
\nu^{1/3}\left(\frac{1-z^2}{4\zeta}\right)^{1/4}
J_{\nu}(\nu z)=\mathcal{A}i(\nu,z)+\epsilon(\nu,z),
\end{equation}
where for $1 \leq z < \infty$ ($-\infty < \zeta \leq 0$), and $1 \leq \nu < \infty$
\begin{equation}
\label{eq1.41}
|\epsilon(\nu,z)|
< \frac{2.0178882463\,\hat{z}_{3}(z)
M(\nu^{2/3}\zeta)}{\nu^{7}z}.
\end{equation}
\end{theorem}

\begin{remark}
\label{remark:Mincreasing}
We have modified the definition of the so-called Airy modulus function $M(x)$ \cite[Chap. 11, Eq. (2.05)]{Olver:1997:ASF} to include all real values of the argument (thereby discarding Eq. (2.04) of that reference). This will be required in our proof. In \cite[Chap. 11, Lemma 5.1]{Olver:1997:ASF} Olver proves that $M(x)$ is increasing for the negative values of $x$ of his definition, but his proof can readily be extended to all $x$ when only using (\ref{eq1.39}). Although our $M(x)$ is unbounded as $x \to \infty$ we shall only use it in the proof for $x$ less than or equal to the small positive value $\eta(1,1) = 0.01711\cdots$ of (\ref{eq1.32}).
\end{remark}

\begin{remark}
\label{remark:AiryErrorSmall}
From \cite[Chap. 11, Eq. (2.07)]{Olver:1997:ASF} $M(x)=\mathcal{O}(|x|^{-1/4})$ as $x \to -\infty$, and from (\ref{eq1.3c}) and (\ref{eq1.19}) $|\zeta| =\mathcal{O}(z^{2/3})$, $\hat{z}_{3}(z)=\mathcal{O}(z^{-5})$ as $z \to \infty$. Thus from (\ref{eq1.41}) $\epsilon(\nu,z) = \mathcal{O}(\nu^{-7}z^{-37/6})$ uniformly for $1 \leq z < \infty$ ($-\infty < \zeta \leq 0$), and $\epsilon(\nu,z) = \mathcal{O}(\nu^{-43/6}z^{-37/6})$ uniformly for $1+\delta \leq z < \infty$ ($\delta > 0$).
\end{remark}

Next for $m=1,2,3,\ldots$ and $\nu \geq1$ define $\hat{z}_{\nu,m}$ implicitly by
\begin{equation}
\label{eq1.42}
\mathcal{Z}_{3}\left(\nu,\hat{z}_{\nu,m}\right)
= \frac{\mathrm{a}_{m}}{\nu^{2/3}},
\end{equation}
and hence $z=\hat{z}_{\nu,m}$ is a zero of $\mathrm{Ai}(\nu^{2/3}\mathcal{Z}_{3}(\nu,z))$.

Let us now define the main error terms that will be bounded. Firstly, on recalling $J_{\nu}(j_{\nu,m})=\mathrm{Ai}(\mathrm{a}_{m})=0$ and referring to (\ref{eq1.38}), (\ref{eq1.40}) and (\ref{eq1.42}), we expect that $\nu^{-1}j_{\nu,m} \approx \hat{z}_{\nu,m}$. Thus let $\hat{\epsilon}_{\nu,m}$ be defined by
\begin{equation}
\label{eq1.45}
\nu^{-1}j_{\nu,m}=\hat{z}_{\nu,m}
+\hat{\epsilon}_{\nu,m}.
\end{equation}
Further, from (\ref{eq1.12}), (\ref{eq1.42}) and \cref{thm:AiryError}, we also expect $\hat{z}_{\nu,m} \approx \mathbf{z}_{\nu,m}$. This is verified by the following, which is proven in \cref{sec4}.

\begin{theorem}
\label{thm:e_bold}
For $m=1,2,3,\ldots$ and $1 \leq \nu < \infty$ 
\begin{equation}
\label{eq1.46}
\hat{z}_{\nu,m} = \mathbf{z}_{\nu,m}
+\mathbf{e}_{\nu,m},
\end{equation}
where $\mathbf{z}_{\nu,m}$ is defined by (\ref{eq1.12}), and the error term satisfies the bounds
\begin{equation}
\label{eq1.48}
 \frac{0.969746 \, z_{m,3}}{\nu^{6}}
<\mathbf{e}_{\nu,m}
< \frac{1.013023 \, z_{m,3}}{\nu^{6}}.
\end{equation}
\end{theorem}

Now, from (\ref{eq1.45}) and (\ref{eq1.46}),
\begin{equation}
\label{eq1.46b}
\nu^{-1}j_{\nu,m} - \mathbf{z}_{\nu,m}
=\mathbf{e}_{\nu,m}+\hat{\epsilon}_{\nu,m},
\end{equation}
and so, on account of (\ref{eq1.48}), in our proof of \cref{thm:main} given in the next section it remains to bound $\hat{\epsilon}_{\nu,m}$. To do so we use the following theorem due to Hethcote \cite{Hethcote:1970:EBA} (with a proof given in his doctoral thesis \cite{Hethcote:1968:PHD}).

\begin{theorem}
\label{thm:hethcote}
In the interval $[b-\rho_{1},b+\rho_{2}]$, suppose $f(z)=g(z)+\epsilon(z)$, where $f(z)$ is continuous, $g(z)$ is differentiable, $g(b)=0$, $\mu=\mathrm{min}|g'(z)|>0$ and
\begin{equation}
\label{eq1.49}
E=\max|\epsilon(z)|
< \min\{|g(b-\rho_{1})|,|g(b+\rho_{2})|\}.
\end{equation}
Then $f(z)$ has a zero $c \in (b-\rho_{1},b+\rho_{2})$ such that $|c-b| \leq E/\mu$.
\end{theorem}

\begin{remark}
In his proof Hethcote assumes an interval with $b$ at the centre, i.e. $\rho_{1}=\rho_{2}$, but his proof can readily be modified for these being unequal as in our application below.
\end{remark}

Based on \cref{thm:AiryError} we shall apply \cref{thm:hethcote} with $b=\hat{z}_{\nu,m}$, $c=\nu^{-1}j_{\nu,m}$ (so that from (\ref{eq1.45}) $c-b=\hat{\epsilon}_{\nu,m}$), and
\begin{equation}
\label{eq1.59}
\begin{split}
f(z) & = \nu^{1/3}\left(1+\dot{\eta}\right)^{1/2}
\left(\frac{1-z^2}{4\zeta}\right)^{1/4}
J_{\nu}(\nu z),
\\
g(z) & = \mathrm{Ai}\left( \nu^{2/3}
\mathcal{Z}_{3}(\nu,z)\right), \quad
\epsilon(z)=
\left(1+\dot{\eta}\right)^{1/2}\epsilon(\nu,z).
\end{split}
\end{equation}

In order to be able to do this we shall utilise certain bounds related to Airy functions, given in the next theorem. In this, (\ref{eq1.53}), (\ref{eq1.54}) and (\ref{eq1.56}), the latter for $m=1,2$, were confirmed by us via explicit computation, with the other more general results proven by Qu and Wong \cite{Qu:1999:BPU} (our notation differs slightly from theirs).

\begin{theorem}
\label{thm:QuWong}
For $m=1,2,3,\ldots$ define $\mathrm{a}_{m,0}$ by (\ref{eq1.50}), and $r_{m}^{-}$ and $r_{m}^{+}$ by
\begin{equation}
\label{eq1.51}
\mathrm{a}_{m}-r_{m}^{-}
=\mathrm{a}_{m,0}\left(
1+\frac{0.01+0.03\,\delta_{m}}{4m-1}\right),
\end{equation}
and
\begin{equation}
\label{eq1.52}
\mathrm{a}_{m}+r_{m}^{+}
=\mathrm{a}_{m,0}\left(
1-\frac{0.01}{4m-1}\right),
\end{equation}
where $\delta_{m}=1$ for $m=1,2$ and is zero otherwise. Then $r_{m}^{-}$ and $r_{m}^{+}$ ($m=1,2,3,\ldots$) are positive,
\begin{equation}
\label{eq1.53}
\min\left\{
\left|\mathrm{Ai}\left(
\mathrm{a}_{1}-r_{1}^{-}\right)\right|,
\left|\mathrm{Ai}\left(
\mathrm{a}_{1}+r_{1}^{+}\right)\right|
\right\}
>9.171267504 \times 10^{-3},
\end{equation}
\begin{equation}
\label{eq1.54}
\min\left\{
\left|\mathrm{Ai}\left(
\mathrm{a}_{2}-r_{2}^{-}\right)\right|,
\left|\mathrm{Ai}\left(
\mathrm{a}_{2}+r_{2}^{+}\right)\right|
\right\}
>9.612776459 \times 10^{-3},
\end{equation}
and for $m=3,4,5,\ldots$
\begin{equation}
\label{eq1.55}
\min\left\{
\left|\mathrm{Ai}\left(
\mathrm{a}_{m}-r_{m}^{-}\right)\right|,
\left|\mathrm{Ai}\left(
\mathrm{a}_{m}+r_{m}^{+}\right)\right|
\right\}
>\frac{3.230051079\times 10^{-3}}
{\sqrt{\pi}\left|\mathrm{a}_{m,0}\right|^{1/4}}.
\end{equation}
Moreover, for $m=1,2,3,\ldots$
\begin{equation}
\label{eq1.56}
\min_{x \in [\mathrm{a}_{m}-r_{m}^{-},
\mathrm{a}_{m}+r_{m}^{+}]}
\left\{
\left|\mathrm{Ai}'(x)\right|
\right\}
>\frac{0.987836345}
{\sqrt{\pi}}
\left|\mathrm{a}_{m,0}\right|^{1/4}.
\end{equation}
\end{theorem}

\begin{remark}
$\mathrm{a}_{m} = \mathrm{a}_{m,0}\{1+\mathcal{O}(m^{-2})\}$ as $m \to \infty$ (see, for example, \cite[Eqs. 9.9.6 and 9.9.18]{NIST:DLMF}), and as such $r_{m}^{\pm}=|\mathrm{a}_{m}|\mathcal{O}(m^{-1})$. 
\end{remark}

Now, for $m=1,2,3,\ldots$, we define $z=w_{\nu,m}^{\pm}$ to be the values that correspond to the end points of the interval in \cref{thm:QuWong}, namely (\ref{eq1.51}) and (\ref{eq1.52}), for the argument $\nu^{2/3}\mathcal{Z}_{3}(\nu,z)$ of the Airy function $g(z)$ in (\ref{eq1.59}). Thus $w_{\nu,m}^{\pm}$ are given implicitly by 
\begin{equation}
\label{eq1.57}
\nu^{2/3}\mathcal{Z}_{3}(\nu,w_{\nu,m}^{+})
=\mathrm{a}_{m} -  r_{m}^{-}, \quad
\nu^{2/3}\mathcal{Z}_{3}(\nu,w_{\nu,m}^{-})
=\mathrm{a}_{m} + r_{m}^{+}.
\end{equation}
With these assigned, and recalling $\hat{z}_{\nu,m}$ is a zero of $g(z)=\mathrm{Ai}(\nu^{2/3}\mathcal{Z}_{3}(\nu,z))$, we define $\rho_{1,2}$ of \cref{thm:hethcote} implicitly in terms of $r_{m}^{\pm}$ via $w_{\nu,m}^{\pm}$ by
\begin{equation}
\label{eq1.58}
w_{\nu,m}^{-}=\hat{z}_{\nu,m} - \rho_{1},
\quad
w_{\nu,m}^{+}=\hat{z}_{\nu,m} + \rho_{2}.
\end{equation}

The advantage of $r_m^{\pm}$ being defined explicitly is that we can directly apply the Airy function estimates from \cref{thm:QuWong} to the corresponding functions in \cref{thm:hethcote}. The price we pay is that the $z$ end points of the interval in \cref{thm:hethcote}, namely $w_{\nu,m}^{\pm}$, are not explicitly given. However, 
for large $\nu$ and/or $m$ we expect $w_{\nu,m}^{\pm}$ to be ``close'' to $z_{m,0}$. In order to obtain simple and explicit error bounds that reflect this we need to estimate $w_{\nu,m}^{\pm}$ in terms of $z_{m,0}$ (which is given by (\ref{eq1.4})). The following provides the required bounds.

\begin{lemma}
\label{lem:wnumBounds}
For $m=1,2,3,\ldots$ and $1 \leq \nu <\infty$
\begin{equation}
\label{eq1.60}
\max\{0.9835\,z_{m,0},1\}
< w_{\nu,m}^{-}< w_{\nu,m}^{+}<1.01354\,z_{m,0}.
\end{equation}
\end{lemma}

\section{Proof of \cref{thm:main}}
\label{sec2}

The following results will be used.

\begin{lemma}
\label{lem:sigma}
For $1<z<\infty$ ($-\infty<\zeta <0$) $\sigma(z)$ is strictly decreasing and $z\sigma(z)$ is strictly increasing.
\end{lemma}

\begin{lemma}
\label{lem:etaBound}
For $-\infty < \zeta \leq 0$  and $1 \leq \nu < \infty$
\begin{equation}
\label{eq1.37}
1 < \left(1+\dot{\eta}\right)^{1/2}
\leq \left(\tfrac{2611707229667}
{2590330640625}\right)^{1/2}
=1.0041177499 \cdots.
\end{equation}
\end{lemma}

\begin{lemma}
\label{lem:MBound}
For $m=1,2,3,\ldots$, $1 \leq \nu < \infty$ and $w_{\nu,m}^{-} \leq z < \infty$
\begin{equation}
\label{eq1.62}
\hat{z}_{3}(z) < (0.9835)^{-5}\,  z_{m,3}
< 1.0867463213\, z_{m,3},
\end{equation}
and
\begin{equation}
\label{eq1.63}
M(\nu^{2/3}\zeta)
< \frac{1.0000277287}{\sqrt{\pi}
\left|\mathrm{a}_{m,0}\right|^{1/4}}.
\end{equation}
\end{lemma}

Now, from our definition (\ref{eq1.59}) of $g(z)$ one has from (\ref{eq1.35a}), (\ref{eq1.22}) and (\ref{eq1.23})
\begin{equation}
\label{eq1.64}
g'(z)= \nu^{2/3}
\left\{-\frac{1}{z\sigma}+\eta'(\nu,z)\right\}
\mathrm{Ai}'\left( \nu^{2/3}
\mathcal{Z}_{3}(\nu,z)\right).
\end{equation}
Next from (\ref{eq1.35}), (\ref{eq1.57}) and (\ref{eq1.58})
\begin{equation*}
w_{\nu,m}^{-} \leq z \leq w_{\nu,m}^{+}
\implies
\mathrm{a}_{m}-r_{m}^{-}\leq \nu^{2/3}\mathcal{Z}_{3}(\nu,z) \leq \mathrm{a}_{m}+r_{m}^{+},
\end{equation*}
and so from (\ref{eq1.34}), (\ref{eq1.56}), (\ref{eq1.60}), (\ref{eq1.64}) and \cref{lem:sigma}
\begin{equation}
\label{eq1.65}
\mu =\min_{z \in [w_{\nu,m}^{-},w_{\nu,m}^{+}]}|g'(z)|
>  \frac{0.9746397231 \, \nu^{2/3}
\left|\mathrm{a}_{m,0}\right|^{1/4}}
{\sqrt{\pi}\,z_{m,0} \,
\sigma(1.01354\,z_{m,0})},
\end{equation}
where we used $0.987836345/1.01354 >0.9746397231$. Next, with the notation $E$ of (\ref{eq1.49}) and $\epsilon(z)$ of (\ref{eq1.59}), we have in this interval, using \cref{thm:AiryError} and \cref{lem:etaBound,lem:wnumBounds,lem:MBound},
\begin{equation}
\begin{split}
\label{eq1.66}
E & =\sup_{z \in [w_{\nu,m}^{-},w_{\nu,m}^{+}]} 
\left\{\left(1+\dot{\eta}\right)^{1/2}
|\epsilon(\nu,z)| \right\}
<1.0041177500
\sup_{z \in [w_{\nu,m}^{-},w_{\nu,m}^{+}]} 
|\epsilon(\nu,z)|
\\ &
<  \frac{2.0262535895}{\sqrt{\pi}\, \nu^{7}
\left|\mathrm{a}_{m,0}\right|^{1/4}}
\sup_{z \in [w_{\nu,m}^{-},w_{\nu,m}^{+}]}
\left\{\frac{\hat{z}_{3}(z)}{z}
\right\}
< \frac{2.0262535895\, z_{m,3}}
{(0.9835)^{6}\,\sqrt{\pi}\, \nu^{7}
\left|\mathrm{a}_{m,0}\right|^{1/4} z_{m,0}}
\\ &
< \frac{2.2389665829\, z_{m,3}}
{\sqrt{\pi}\, \nu^{7}
\left|\mathrm{a}_{m,0}\right|^{1/4} z_{m,0}}.
\end{split}
\end{equation}

For the criteria (\ref{eq1.49}) of \cref{thm:hethcote} we have from (\ref{eq1.50}), (\ref{eq1.16}) and the final bound of (\ref{eq1.66})
\begin{equation*}
E < \frac{2.2389665829\, \hat{z}_{3}(1)}
{\sqrt{\pi}\,
\left|\mathrm{a}_{1,0}\right|^{1/4}}
= 0.00045 \cdots
\quad (m=1,2,3,\ldots,\,1 \leq  \nu < \infty),
\end{equation*}
where we have simplified by using the inequalities $z_{m,0} > 1$ and $z_{m,3}=\hat{z}_{3}(z_{m,0}) < \hat{z}_{3}(1)$ (see \cref{lem:zmjMonotonic}). This is smaller than the lower bounds of (\ref{eq1.53}) and (\ref{eq1.54}), and hence for $m=1,2$ the requirement is met. 

To verify the same is true for other values of $m$ we similarly obtain
\begin{equation*}
E < \frac{2.2389665829\, \hat{z}_{3}(1)}
{\sqrt{\pi}\,
\left|\mathrm{a}_{m,0}\right|^{1/4}}
= \frac{0.00099 \cdots}{\sqrt{\pi}\,
\left|\mathrm{a}_{m,0}\right|^{1/4}}
\quad (m=1,2,3,\ldots,\,1 \leq  \nu < \infty).
\end{equation*}
This proves that $E$ is smaller than the lower bound of (\ref{eq1.55}), and hence for $m=3,4,5,\ldots$ the requirement (\ref{eq1.49}) is again met.

Returning to the general case, again from the last bound in (\ref{eq1.66}) and dividing by (\ref{eq1.65}), we deduce from \cref{thm:hethcote} (recalling $c-b=\hat{\epsilon}_{\nu,m}$) that the following holds
\begin{equation}
\label{eq1.71}
|\hat{\epsilon}_{\nu,m}|
\leq \frac{E}{\mu} 
<  \frac{z_{m,3} \, \chi_{m}}{\nu^{23/3}}
\quad (m=1,2,3,\ldots,\,1 \leq  \nu < \infty),
\end{equation}
where $\chi_{m}$ is given by (\ref{eq1.71a}).

Finally, consider $\nu^{-5/3}\chi_{m}$ appearing in (\ref{eq1.72}). Taking into consideration (\ref{eq1.50}) and that $\sigma(z)$ is decreasing for $1 < z < \infty$ (\cref{lem:sigma}), and $z_{m,0}$ is increasing as a function of $m=1,2,3,\ldots$ (see (\ref{eq1.4})), the largest value of $\nu^{-5/3}\chi_{m}$ for each fixed $\nu \geq 1$ is when $m=1$. For this value we can regard it as a function of the single variable $\nu$, bearing in mind from (\ref{eq1.4}) that $z_{1,0}$ also depends on $\nu$. Numerically we confirm $\chi_{1}$ is a decreasing function for $1 \leq \nu < \infty$, and so with $\tilde{z}:=1.01354 \,z_{1,0}=3.87444\cdots$ when $\nu=1$ (evaluated from (\ref{eq1.4})) we deduce that
\begin{equation}
\label{eq1.75a}
\frac{\chi_{m}}{\nu^{5/3}}
\leq \frac{2.2972248409 \,
\sigma(\tilde{z})}
{\left|\mathrm{a}_{1,0}\right|^{1/2}}
= 0.62034 \cdots
\quad (m=1,2,3,\ldots,\,1 \leq  \nu < \infty).
\end{equation}

In conclusion, from (\ref{eq1.48}), (\ref{eq1.46b}) and (\ref{eq1.71}) we obtain (\ref{eq1.72}), with (\ref{eq1.75a}) establishing that the lower bound is indeed positive. This completes the proof.

\section{Proof of \cref{thm:AiryError}}
\label{sec3}
We begin with the differential equation
\begin{equation}
\label{eq3.1}
d^{2}W/d \zeta^{2}
=\left\{\nu^{2}\zeta+\psi(\zeta)\right\}W,
\end{equation}
where
\begin{equation}
\label{eq3.2}
\psi(\zeta)=\frac{5}{16\zeta^{2}}
+\frac{\zeta z^{2}(z^{2}+4)}
{4\left(z^{2}-1\right)^{3}}.
\end{equation}
As shown in \cite[Chap. 11, Sect. 10]{Olver:1997:ASF} this has solutions $W(\nu,\zeta)=\{(1-z^2)/\zeta\}^{1/4}\mathscr{C}_{\nu}(\nu z)$, where $\mathscr{C}_{\nu}(z)$ is a solution of Bessel's equation (\cite[Eq. 10.2.1]{NIST:DLMF}). The equation (\ref{eq3.1}) is precisely of the form for which \cite[Chap. 11, Thm. 9.1]{Olver:1997:ASF} is applicable, providing uniform asymptotic expansions involving Airy functions and their derivatives.

In place of Olver's expansions we assume a solution of the form
\begin{equation}
\label{eq3.3}
W(\nu, z)
=\mathcal{A}i(\nu,z)
+\epsilon(\nu,z),
\end{equation}
where $\mathcal{A}i(\nu,z)$ is defined by (\ref{eq1.38}). Our goal is to bound the error term $\epsilon(\nu,z)$ uniformly for $-\infty < \zeta \leq 0$ (corresponding to $1 \leq z < \infty$).

Recalling that dots represent differentiation with respect to $\zeta$ we obtain on inserting (\ref{eq3.3}) into (\ref{eq3.1}), and referring to Airy's equation (\cite[Eq. 9.2.1]{NIST:DLMF}),
\begin{equation}
\label{eq3.6}
\ddot{e}(\nu,\zeta)-\nu^{2}\zeta e(\nu,\zeta)
=\psi(\zeta) e(\nu,\zeta) +\gamma(\nu,\zeta)
\mathscr{A}(\nu,\zeta),
\end{equation}
where $e(\nu,\zeta)=\epsilon(\nu,z(\zeta))$, $\mathscr{A}(\nu,\zeta)=\mathcal{A}i(\nu,z(\zeta))$ and
\begin{equation}
\label{eq3.7}
\gamma=\psi-\nu^{2}\left\{\eta+(\zeta+\eta)\left(2+\dot{\eta}\right)\dot{\eta}\right\}
-\frac{3 {\ddot{\eta}}^{2}-2(1+\dot{\eta})
\dddot{\eta}}
{4\left(1+\dot{\eta}\right)^{2}}.
\end{equation}
From (\ref{eq1.3c}), (\ref{eq1.23}) and (\ref{eq3.2})
\begin{equation}
\label{eq3.33}
\gamma(\nu,\zeta)
=\frac{7672012\, (12^{2/3})}{2657205 \, z^{22/3}}
\left\{1+\mathcal{O}\left(\frac{1}{z}\right)
\right\}
\quad (z \to \infty).
\end{equation}
We later shall show that  $\gamma(\nu,\zeta)=\mathcal{O}(\nu^{-6})$ as $\nu \to \infty$, uniformly for $-\infty < \zeta \leq 0$.

Now on solving (\ref{eq3.6}) by variation of parameters we obtain in the standard manner the integral equation
\begin{equation}
\label{eq3.8}
e(\nu,\zeta)=\frac{\pi}{\nu^{2/3}}
\int_{-\infty}^{\zeta}K(\zeta,t)
\left\{\gamma(\nu,t)\mathscr{A}(\nu,t)
+\psi(t)e(\nu,t)\right\} dt,
\end{equation}
where
\begin{equation*}
K(\zeta,t)=
\mathrm{Ai}(\nu^{2/3}\zeta)
\mathrm{Bi}(\nu^{2/3}t)
-\mathrm{Ai}(\nu^{2/3}t)
\mathrm{Bi}(\nu^{2/3}\zeta).
\end{equation*}

For $-\infty< t \leq \zeta \leq 0$ we have from \cite[Chap. 11, Eq. (3.14)]{Olver:1997:ASF}, noting that in this equation the so-called weight functions $E$ are identically equal to $1$ in the oscillatory intervals we are considering,
\begin{equation}
\label{eq3.10}
|K(\zeta,t)| \leq
M(\nu^{2/3}\zeta) M(\nu^{2/3}t),
\end{equation}
where $M(x)$ is defined by (\ref{eq1.39}) (see also \cref{remark:Mincreasing}). We note that from \cite[Chap. 11, Eq. (2.07)]{Olver:1997:ASF}
\begin{equation}
\label{eq3.11}
M(x) \sim \pi^{-1/2}|x|^{-1/4}
\quad  (x \to -\infty).
\end{equation}

Identifying (\ref{eq3.8}) with \cite[Chap. 6, Eq. (10.01)]{Olver:1997:ASF} we replace Olver's $\xi$ by $\zeta$, set $h(\zeta)=e(\nu,\zeta)$ and, taking into account (\ref{eq3.10}),
\begin{equation}
\label{eq3.12}
\begin{split}
\mathsf{K}(\zeta,t) & = \pi\nu^{-2/3}|t|^{1/2}K(\zeta,t), \;
P_{0}(\zeta)=\pi \nu^{-1}M(\nu^{2/3}\zeta), \;
\\
Q(t) & = \nu^{1/3}|t|^{1/2}M(\nu^{2/3}t), \;
\phi(t)=|t|^{-1/2}\gamma(\nu,t), \;
\\
J(t) & = \mathscr{A}(\nu,t), \;
\psi_{0}(t)=|t|^{-1/2}\psi(t), \;
\psi_{1}(t)=0.
\end{split}
\end{equation}
As a result from \cite[Chap. 6, Thm. 10.1]{Olver:1997:ASF} we obtain the following bound.
\begin{lemma}
For $-\infty < \zeta \leq 0$ and $\nu >0$
\begin{equation}
\label{eq3.13}
|e(\nu,\zeta)|
\leq
\kappa \nu^{-1}M(\nu^{2/3}\zeta)
\Phi(\nu,\zeta)
\exp \left\{\nu^{-1} \Psi_{0}(\zeta)\right\},
\end{equation}
where
\begin{equation}
\label{eq3.14}
\kappa
=\sup_{-\infty<\zeta \leq 0}
\left\{\pi  \nu^{1/3}|\zeta|^{1/2}M(\nu^{2/3}\zeta)
\left|\mathscr{A}(\nu,\zeta)\right|\right\},
\end{equation}
\begin{equation}
\label{eq3.16}
\Phi(\nu,\zeta)
= \int_{-\infty}^{\zeta}
|t^{-1/2}\gamma(\nu,t)| dt,
\end{equation}
and
\begin{equation}
\label{eq3.17}
\Psi_{0}(\zeta)
= \int_{-\infty}^{\zeta}
|t^{-1/2}\psi(t)| dt.
\end{equation}
\end{lemma}

\begin{remark}
In the notation of \cite[Chap. 6, Thm. 10.1]{Olver:1997:ASF} we computed
\begin{equation*}
\label{eq3.15}
\kappa_{0}
= \sup_{-\infty<x \leq 0}
\left\{ \pi  |x|^{1/2}
M^{2}(x)\right\} =1,
\end{equation*}
with this supremum attained at $x = - \infty$ (see (\ref{eq3.11})).
\end{remark}

Let us now simplify various terms in the bound (\ref{eq3.13}). Firstly from (\ref{eq1.3})
\begin{equation}
\label{eq3.19}
\frac{d \zeta}{|\zeta|^{1/2}}
=-\frac{\left(z^{2}-1\right)^{1/2}
dz}{z |\zeta|},
\end{equation}
and hence again from (\ref{eq1.3}), along with (\ref{eq3.2}) and (\ref{eq3.17}) and numerical integration,
\begin{multline}
\label{eq3.20}
\Psi_{0}(\zeta)
\leq \Psi_{0}(0)
= \int_{1}^{\infty}
\frac{\left(z^{2}-1\right)^{1/2}
|\psi(\zeta(z))| dz}{z |\zeta|}
\\
= 0.0434514175\cdots
\quad  (-\infty < \zeta \leq 0).
\end{multline}
With this, one obtains for the exponential term in (\ref{eq3.13}), the simplified bound
\begin{multline}
\label{eq3.21}
\exp \left\{\nu^{-1} 
\Psi_{0}(\zeta)\right\}
\leq \exp \left\{\Psi_{0}(0)\right\}
= 1.0444092531 \cdots
\\
\quad (-\infty < \zeta \leq 0, \, 1 \leq \nu <\infty).
\end{multline}

Next, from (\ref{eq1.38}), (\ref{eq3.14}) and \cref{lem:etaBound} we deduce that $\kappa \leq \kappa_{1}\kappa_{2}$, where
\begin{equation}
\label{eq3.24}
\kappa_{1}
=\sup_{-\infty<\zeta \leq 0} \left\{
\left|1+\dot{\eta}\right|^{-1/2} \right\}<1
\quad  (1 \leq \nu < \infty),
\end{equation}
and
\begin{equation}
\label{eq3.25}
\kappa_{2}
=\sup_{-\infty<\zeta \leq 0}
\left\{\pi  \nu^{1/3}|\zeta|^{1/2}
M(\nu^{2/3}\zeta)
\left|\mathrm{Ai}\left(\nu^{2/3}\mathcal{Z}_{3}\right)
\right|\right\}.
\end{equation}
Now from (\ref{eq1.22}), (\ref{eq1.23}), (\ref{eq1.33}) and (\ref{eq1.39}), for $-\infty < \zeta  \leq 0$ and $1 \leq \nu < \infty$,
\begin{equation*}
\label{eq3.26}
\left|\mathrm{Ai}\left(\nu^{2/3}\mathcal{Z}_{3}\right)
\right| \leq
M\left(\nu^{2/3}\mathcal{Z}_{3}\right),
\quad
\nu^{2/3}\mathcal{Z}_{3}
\leq \nu^{2/3}\zeta+\eta(1,1)
=\nu^{2/3}\zeta
+\tfrac {44873962351}{3302530481250} 2^{1/3}.
\end{equation*}
Consequently, with the fact that $M(x)$ is monotonically increasing for all $x$ (see \cref{remark:Mincreasing}),
\begin{equation}
\label{eq3.28}
\kappa_{2}
\leq \sup_{-\infty<x \leq 0}
\left\{ \pi  |x|^{1/2} M(x)
M\left(x+\eta(1,1)\right)\right\}
=1.000273093257\cdots,
\end{equation}
with the supremum attained at $x=-10.44187\cdots$. Thus, on recalling $\kappa \leq \kappa_{1}\kappa_{2}$, from (\ref{eq3.21}), (\ref{eq3.24}) and (\ref{eq3.28}) we arrive at
\begin{equation}
\label{eq3.55}
\kappa \exp \left\{\nu^{-1}
\Psi_{0}(\zeta)\right\}
< 1.0446944743
\quad (-\infty < \zeta \leq 0,\, 1 \leq \nu <\infty).
\end{equation}

It remains in (\ref{eq3.13}) to obtain a simple and computable bound for the integral $\Phi$ given by (\ref{eq3.16}), recalling that here $\gamma$ is defined by (\ref{eq3.7}). Although this function only involves explicit elementary functions, as we see shortly it is very unwieldy containing many terms, particularly when the higher $\zeta$ derivatives of $\eta$ are converted to $z$ derivatives via (\ref{eq1.35a}) and the chain rule. In obtaining our simplified bound we also shall confirm that $\gamma(\nu,\zeta)=\mathcal{O}(\nu^{-6})$ as $\nu \to \infty$ uniformly for $-\infty < \zeta \leq 0$ ($1 \leq z < \infty$), with the same then being true of $\Phi(\nu,\zeta)$.

We begin by using for $\dot{\eta} \neq -1$ the simple algebraic identity
\begin{equation*}
\frac{1}{(1+\dot{\eta})^{2}}
=1-2\dot{\eta}
+\frac{\dot{\eta}^{2}(3+2\dot{\eta})}
{(1+\dot{\eta})^{2}},
\end{equation*}
and as a result from (\ref{eq3.7}) we recast $\gamma=\gamma_{1}+\gamma_{2}$ where
\begin{equation}
\label{eq3.34}
\gamma_{1}=\psi-\nu^{2}\left\{\eta+(\zeta+\eta)\left(2+\dot{\eta}\right)\dot{\eta}\right\}
-\tfrac{1}{4} (1-2\dot{\eta})
\{3\ddot{\eta}^{2}-2(1+\dot{\eta})\dddot{\eta}\},
\end{equation}
and
\begin{equation}
\label{eq3.35}
\gamma_{2}=
-\frac{\dot{\eta}^{2}
(3+2\dot{\eta})
\{3{\ddot{\eta}}^{2}-2(1+\dot{\eta})
\dddot{\eta}\}}
{4\left(1+\dot{\eta}\right)^{2}}.
\end{equation}

From using (\ref{eq1.23}) - (\ref{eq1.26}) and (\ref{eq3.2}) we find that the $\mathcal{O}(\nu^{-2j})$ ($j=0,1,2$) terms cancel in $\gamma_{1}$. Thus $\gamma_{1}=\mathcal{O}(\nu^{-6})$ as $\nu \to \infty$, and in fact
\begin{equation}
\label{eq3.36}
\gamma_{1} = \frac{1}{\nu^{6}}
\sum_{j=0}^{6}\frac{q_{j}}{\nu^{2j}},
\end{equation}
where from (\ref{eq1.23}) and (\ref{eq3.34}) the functions $q_{j}$ are independent of $\nu$ and can be explicitly expressed in terms of $\zeta$, and for $l=1,2,3$, $\Upsilon_{l}$, $\dot{\Upsilon}_{l}$, $\ddot{\Upsilon_{l}}$ and $\dddot{\Upsilon}_{l}$. The leading term is given by
\begin{equation*}
\begin{split}
q_{1}= &
-\left( \dddot{{\Upsilon}}_{1}+\Upsilon_{2} \right) {{ \dot{\Upsilon}_{1}}}^{2}-
 \tfrac{1}{2}\left( 4 \zeta{ \dot{\Upsilon}_{3}}
 +4\Upsilon_{1}{ \dot{\Upsilon}_{2}}
-3 \ddot{{\Upsilon}}_{1}^{2}+4\Upsilon_{3}+{ \dddot{{\Upsilon}}_{2}} \right) 
{ \dot{\Upsilon}_{1}}
-\zeta{ \dot{\Upsilon}_{2}}^{2}
\\ &
- \tfrac{1}{2}\left(\dddot{{\Upsilon}}_{1}
+2\Upsilon_{2} \right) 
\dot{\Upsilon}_{2}-2\Upsilon_{1}{ \dot{\Upsilon}_{3}}
-\tfrac{3}{2}\ddot{\Upsilon}_{1}\ddot{{\Upsilon}}_{2}
+\tfrac{1}{2}\dddot{\Upsilon}_{3},
\end{split}
\end{equation*}
with the others explicitly obtainable but not recorded here. As we noted above, all the $\zeta$ derivatives can be converted to $z$ derivatives via (\ref{eq1.35a}) and the chain rule. For example, for the first derivatives, we have $\dot{\Upsilon}_{l} = -z \sigma \Upsilon_{l}'$, with as usual dots and primes denoting $\zeta$ and $z$ derivatives, respectively.

Now consider $\gamma_{2}$ given by (\ref{eq3.35}). We do not have to expand in inverse powers of $\nu$ like we did for $\gamma_{1}$ in (\ref{eq3.36}) since we can see directly from (\ref{eq1.23}) that it is also $\mathcal{O}(\nu^{-6})$ as $\nu \to \infty$, since there are no cancellations of lower order terms. For $1 \leq \nu < \infty$ we can bound this function using (\ref{eq1.23}) and the simple inequalities $|\dot{\eta}| \leq \nu^{-2} g_{1}$, $|\ddot{\eta}| \leq \nu^{-2} g_{2}$ and $|\dddot{\eta}| \leq \nu^{-2} g_{3}$, where
\begin{equation}
\label{eq3.39}
g_{1}
= |\dot{\Upsilon}_{1}|
+|\dot{\Upsilon}_{2}|
+|\dot{\Upsilon}_{3}|,
\end{equation}
with $g_{2}$ and $g_{3}$ similarly defined with each $\dot{\Upsilon}$ replaced by $\ddot{\Upsilon}$ and $\dddot{\Upsilon}$, respectively.

Thus from (\ref{eq3.35}) and (\ref{eq3.36}), and recalling that $\gamma=\gamma_{1}+\gamma_{2}$, we have
\begin{equation}
\label{eq3.40}
\begin{split}
|\gamma(\nu,\zeta)| 
& =  |\gamma_{1}(\nu,\zeta)+\gamma_{2}(\nu,\zeta)|
\leq |\gamma_{1}(\nu,\zeta)|+|\gamma_{2}(\nu,\zeta)|
\\ &
\leq \nu^{-6}
\left\{G_{1}(\zeta)+G_{2}(\zeta)\right\}
\quad (-\infty < \zeta \leq 0, \, 1 \leq \nu < \infty),
\end{split}
\end{equation}
where
\begin{equation}
\label{eq3.41}
G_{1} = 
\sum_{j=0}^{6}|q_{j}|,
\end{equation}
and
\begin{equation}
\label{eq3.42}
G_{2}=\frac{g_{1}^{2}
(3+2g_{1})\{3g_{2}^{2}+2(1+g_{1})g_{3}\}}
{4\left(1-g_{1}\right)^{2}}.
\end{equation}
For the denominator of $G_{2}$ one can show numerically that $g_{1}$ is less than 1 for $-\infty < \zeta \leq 0$ ($1 \leq z < \infty$), and in fact it attains a maximum value $0.0158\cdots$ at $z=1$. 

Now from (\ref{eq3.16}) and (\ref{eq3.19})
\begin{equation}
\label{eq3.43}
\Phi(\nu,\zeta)
= \int_{z}^{\infty}
\frac{\left(t^{2}-1\right)^{1/2}
|\gamma(\nu,\zeta(t))| dt}{t |\zeta(t)|},
\end{equation}
and hence from (\ref{eq3.40})
\begin{equation}
\label{eq3.44}
\Phi(\nu,\zeta)
\leq \frac{G(z)}{\nu^{6}}
\quad (-\infty < \zeta \leq 0, \,1 \leq \nu < \infty),
\end{equation}
where
\begin{equation}
\label{eq3.45}
G(z) =
\int_{z}^{\infty}
\Lambda(t)dt,  \quad
\Lambda(z)=\frac{\left(z^{2}-1\right)^{1/2}
\left\{G_{1}(\zeta)+G_{2}(\zeta)\right\}}
{z |\zeta|}.
\end{equation}
Note that $\Lambda(z)=\mathcal{O}\{(z-1)^{-1/2}\}$ as $z \to 1$ since from (\ref{eq1.3b}) $\zeta^{-1} = \mathcal{O}\{(z-1)^{-1}\}$.

\begin{figure}[htbp]
 \centering
 \includegraphics[width=0.7\textwidth,
 keepaspectratio]
 {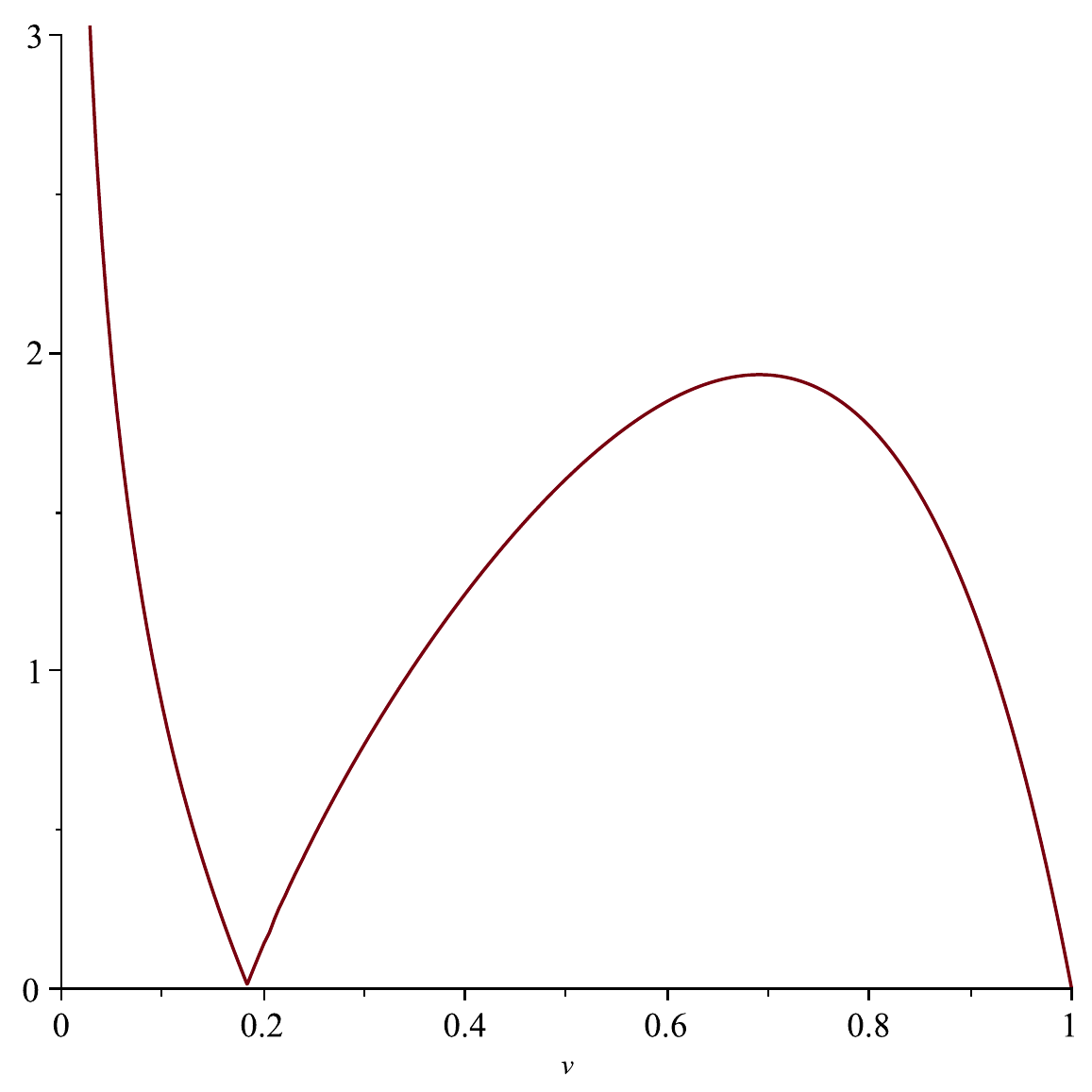}
 \caption{Graph of $p_{1}(v)$ for $0 < v \leq 1$}
 \label{fig:figGoverDeriv}
\end{figure}

Unlike (\ref{eq3.20}), we want our bound on $\Phi(\nu,\zeta)$ to vanish as $z \to \infty$, and in particular to take into account (\ref{eq3.33}) and (\ref{eq3.16}). To this end we shall compare $\nu^{6}\Phi(\nu,\zeta)$ with $z^{-1}\hat{z}_{3}(z)$ for $1 \leq z < \infty$, via (\ref{eq3.44}) and (\ref{eq3.45}). Thus consider
\begin{equation}
\label{eq3.49}
p_{1}(v)
=\frac{\Lambda((1-v)^{-1})}{p_{0}(v)}
\quad   (0 < v <1),
\end{equation}
where
\begin{equation*}
p_{0}(v)
=-\left . \left(
\frac{d}{d z}\frac{\hat{z}_{3}(z)}{z} \right )
\right \vert_{z=(1-v)^{-1}}
\quad   (0 \leq v <1),
\end{equation*}
this latter function being readily shown to be positive for $0 \leq v <1$. The graph of $p_{1}(v)$ is shown in \cref{fig:figGoverDeriv}. This has a vertical asymptote at $v=0$ due to the square root singularity of $\Lambda(z)$ at $z=1$ ($v=0$). It has a local maximum of $1.9315582649\cdots$ at $v = 0.69098\cdots$ (corresponding to $z=3.23605\cdots$). We find that it also attains this value at $v=v_{1}:=0.05151\cdots$ (with the corresponding value $z=x_{1}:=(1-v_{1})^{-1}=1.05430\cdots$). From the graph it is then evident that $p_{1}(v) < 1.9315582650$ for $v_{1} \leq v < 1$. 

Thus replacing $z$ by $t$ and taking into account (\ref{eq3.49}) yields
\begin{equation}
\label{eq3.51}
\Lambda(t) < -1.9315582650 \frac{d}{d t}
\left ( \frac{\hat{z}_{3}(t)}{t} \right )
\quad   (x_{1} \leq t < \infty).
\end{equation}
Then from (\ref{eq3.44}), (\ref{eq3.45}) and (\ref{eq3.51}), one arrives upon integration
\begin{equation}
\label{eq3.52}
\Phi(\nu,\zeta) \leq \frac{1}{\nu^{6}}
\int_{z}^{\infty}\Lambda(t)dt
< \frac{1.9315582650\,\hat{z}_{3}(z)}{\nu^{6}z}
\quad   (x_{1} \leq z < \infty,\, 1 \leq \nu < \infty).
\end{equation}

We can extend this to $1 \leq z < \infty$ as follows. For $1 \leq z \leq x_{1}$, on noting from \cref{lem:zmjMonotonic} that $\hat{z}_{3}(z)$ is decreasing, as is $G(z)$ in this interval (see (\ref{eq3.45})), we deduce from (\ref{eq3.44}) that
\begin{equation}
\label{eq3.53}
\frac{z \Phi(\nu,\zeta)}{\hat{z}_{3}(z)} 
\leq \frac{ x_{1} G(1)}{\nu^{6}\hat{z}_{3}(x_{1})}
\quad   (1 \leq z \leq x_{1},\, 1 \leq \nu < \infty).
\end{equation}

We perform a numerical integration on (\ref{eq3.45}) to compute $G(1)$, aided by the change variable $t=(1-v^2)^{-1}$ ($\implies dt=2v(1-v^2)^{-2} dv$), which not only removes the square root singularity of the integrand at $t=1$ ($v=0$), but also maps $1 \leq t < \infty$ to the finite integration interval $0 \leq v <1$. As a result we find $G(1)=0.00061\cdots$. Using this value in (\ref{eq3.53}), as well as recalling $x_{1}=1.05430\cdots$ and a straightforward calculation from (\ref{eq1.13}) yielding $\hat{z}_{3}(x_{1})=0.00041\cdots$, we arrive at
\begin{equation*}
\Phi(\nu,\zeta)
< \frac{1.543\,\hat{z}_{3}(z)}{\nu^{6}z}
\quad   (1 \leq z \leq x_{1},\, 1 \leq \nu < \infty).
\end{equation*}
Thus the second bound of (\ref{eq3.52}) indeed holds for $1 \leq z < \infty$ ($-\infty < \zeta \leq 0$). Consequently, on combining (\ref{eq3.13}), (\ref{eq3.55}) and (\ref{eq3.52}) (in which $x_{1}$ can now be replaced by $1$), the bound (\ref{eq1.41}) is established.

Finally, we prove the relation (\ref{eq1.40}). Now from (\ref{eq1.3c}), (\ref{eq1.22}), (\ref{eq1.23}), and (\ref{eq1.29}) - (\ref{eq1.31})
\begin{equation}
\label{eq1.31a}
\tfrac{2}{3}|\mathcal{Z}(\nu,z)|^{3/2}
=\tfrac{2}{3}|\zeta|^{3/2}
+\mathcal{O}\left(z^{-1}\right)
=z -\tfrac{1}{2}\pi
+\mathcal{O}\left(z^{-1}\right)
\quad (z \to \infty).
\end{equation}
Then consider the asymptotic solution $W(\nu, z)$ given by (\ref{eq3.3}). For this we have, using (\ref{eq1.31a}) and the asymptotic expansion for the Airy function of large negative argument \cite[Eq. 9.7.9]{NIST:DLMF}, and assuming $\nu >0$,
\begin{equation}
\label{eq3.56a}
W(\nu, z)
=\frac{\cos\left(\nu z-\tfrac{1}{2}\nu \pi 
-\tfrac{1}{4}\pi\right)}
{ \sqrt{\pi} \nu^{1/6}|\zeta|^{1/4} }
+\mathcal{O}\left(\frac{1}{z|\zeta|^{1/4}}\right)
\quad   (\zeta \to -\infty),
\end{equation}
noting that $\epsilon(\nu,z) = \mathcal{O}(z^{-37/6})= \mathcal{O}(\zeta^{-37/4})$ as $\zeta \to -\infty$ (see (\ref{eq1.3c}) and \cref{remark:AiryErrorSmall}).

On the other hand, from the well-known asymptotic behaviour of Bessel functions of large argument (see for example \cite[Eq. 10.17.3]{NIST:DLMF}), we find the function on the LHS of (\ref{eq1.40}) has the identical oscillatory behavior at $\zeta=-\infty$.
The claimed identity then is verified since both of these functions are solutions of the differential equation (\ref{eq3.1}) having this unique behaviour.

\section{Proof of \cref{thm:e_bold}}
\label{sec4}
Firstly we require the following.
\begin{lemma}
\label{lem:zhatbound}
For $m=1,2,3,\ldots$ and $1 \leq \nu < \infty$,  $z=\hat{z}_{\nu,m}$ (defined by (\ref{eq1.42})) is a unique and simple zero  in $(1, \infty)$ of the function
\begin{equation}
\label{eq1.43}
F_{m}(\nu,z)
=\mathcal{Z}_{3}(\nu,z) - \frac{\mathrm{a}_{m}}{\nu^{2/3}},
\end{equation}
and satisfies the bounds
\begin{equation}
\label{eq1.44}
1 < z_{m,0} < \hat{z}_{\nu,m} < z_{m,0} +\tfrac{1}{73}.
\end{equation}
\end{lemma}

Now, on denoting $F_{m}'(\nu,z)= \partial  F_{m}(\nu,z) / \partial z$, we have from Taylor's remainder theorem
\begin{equation}
\label{eq2.2}
0=F_{m}(\nu,\hat{z}_{\nu,m})=
F_{m}(\nu,\mathbf{z}_{\nu,m}+\mathbf{e}_{\nu,m})
=F_{m}(\nu,\mathbf{z}_{\nu,m})
+\mathbf{e}_{\nu,m}F_{m}'(\nu,\tau),
\end{equation}
for some $\tau$ lying between $\mathbf{z}_{\nu,m}$ and $\hat{z}_{\nu,m}$, so that
\begin{equation*}
1 \leq \mathbf{z}_{\nu,m}<\tau<\hat{z}_{\nu,m} \;
\text{ or } \;
1 < \hat{z}_{\nu,m}<\tau<\mathbf{z}_{\nu,m}.
\end{equation*}
It then follows from (\ref{eq1.21}) and (\ref{eq1.44})
\begin{equation}
\label{eq2.4}
z_{m,0} <\tau<z_{m,0}+\tfrac{1}{73}.
\end{equation}

Next, from (\ref{eq2.2}),
\begin{equation}
\label{eq2.5}
\mathbf{e}_{\nu,m}
=-\frac{F_{m}(\nu,\mathbf{z}_{\nu,m})}
{F_{m}'(\nu,\tau)},
\end{equation}
assuming $F_{m}'(\nu,\tau) \neq 0$, which we shall later show to be true (see (\ref{eq2.6}) below). Consider the numerator of this, and let $\delta = \nu^{-2}$. Then from the definition (\ref{eq1.12}) we have  $\mathbf{z}_{\nu,m} = z_{m,0}+\delta z_{m,1}+\delta^{2} z_{m,2}$, and denote
\begin{equation}
\label{eq2.10a}
\breve{\zeta}(\delta)
:=\zeta(\mathbf{z}_{\nu,m})
=\zeta\left(z_{m,0}+\delta z_{m,1}
+\delta^{2} z_{m,2}\right),
\end{equation}
and with a similar notation $\breve{\Upsilon}_{j}(\delta)
:=\Upsilon_{j}(\mathbf{z}_{\nu,m})$ ($j=1,2$). Thus, from (\ref{eq1.22}), (\ref{eq1.23}) and (\ref{eq1.43}),
\begin{equation}
\label{eq2.10b}
F_{m}(\nu,\mathbf{z}_{\nu,m})
=\breve{\zeta}(\delta)
-\nu^{-2/3} \mathrm{a}_{m}
+\delta \breve{\Upsilon}_{1}(\delta)
+\delta^{2} \breve{\Upsilon}_{2}(\delta)
+\delta^{3} \Upsilon_{3}(\mathbf{z}_{\nu,m}).
\end{equation}

Recall from (\ref{eq1.5}) that $\zeta(z_{m,0})=\zeta_{m,0}$, and let $\Upsilon_{j}(z_{m,0})=\Upsilon_{m,j}$ ($j=1,2$), and similarly for the derivatives; for example $\zeta'(z_{m,0})=\zeta_{m,0}'$, $\Upsilon_{j}'(z_{m,0})=\Upsilon_{m,j}'$, etc. Thus from (\ref{eq1.5}) and (\ref{eq2.10a}) we have $\breve{\zeta}(0)=\zeta_{m,0}=\nu^{-2/3} \mathrm{a}_{m}$, and so by Taylor's theorem, assuming $1 \leq \nu < \infty$ ($0 < \delta \leq 1$),
\begin{equation}
\label{eq2.10}
\breve{\zeta}(\delta) -\nu^{-2/3} \mathrm{a}_{m}
=\delta \breve{\zeta}'(0)
+\tfrac{1}{2} \delta^{2} \breve{\zeta}''(0)
+\tfrac{1}{6} \delta^{3} \breve{\zeta}'''(\delta'''),
\end{equation}
where $0< \delta''' < \delta$. Here $\breve{\zeta}'(\delta)=d\breve{\zeta}(\delta)/d \delta$ and similarly for the higher $\delta$ derivatives, as well as $\breve{\Upsilon}_{1}''(\delta)$ and $\breve{\Upsilon}_{2}'(\delta)$ below.

Now, upon explicit differentiation using the chain rule, we have from (\ref{eq2.10a})
\begin{equation*}
\breve{\zeta}'(0)
=z_{m,1} \zeta'_{m,0}, \quad
\breve{\zeta}''(0)
=z_{m,1}^{2} \zeta''_{m,0}
+2 z_{m,2} \zeta'_{m,0},
\end{equation*}
and hence from (\ref{eq2.10})
\begin{equation}
\label{eq2.10c}
\breve{\zeta}(\delta) -\nu^{-2/3} \mathrm{a}_{m}
=\delta z_{m,1} \zeta'_{m,0}
+\tfrac{1}{2} \delta^{2} \left\{
z_{m,1}^{2} \zeta''_{m,0}
+2 z_{m,2} \zeta'_{m,0} \right\}
+\tfrac{1}{6} \delta^{3} \breve{\zeta}'''(\delta''').
\end{equation}

In a similar manner
\begin{equation}
\label{eq2.12}
\delta \breve{\Upsilon}_{1}(\delta)
=\delta \Upsilon_{m,1}+\delta^{2} z_{m,1} \Upsilon_{m,1}'
+\tfrac{1}{2} \delta^{3} \breve{\Upsilon}_{1}''(\delta''),
\end{equation}
and
\begin{equation}
\label{eq2.13}
\delta^{2} \breve{\Upsilon}_{2}(\delta)
=\delta^{2} \Upsilon_{m,2}+\delta^{3}
\breve{\Upsilon}_{2}'(\delta'),
\end{equation}
where $\delta'$ and $\delta''$ are some numbers that lie in $(0,\delta)$. Thus we have for the first four terms on the RHS of (\ref{eq2.10b})
\begin{equation}
\label{eq2.14}
\breve{\zeta}(\delta)-\nu^{-2/3} \mathrm{a}_{m}
+\delta \breve{\Upsilon}_{1}(\delta)
+\delta^{2} \breve{\Upsilon}_{2}(\delta)
=\tfrac{1}{6} \delta^3
\left\{\breve{\zeta}'''(\delta''')
+3\breve{\Upsilon}_{1}''(\delta'')
+6\breve{\Upsilon}_{2}'(\delta')\right\},
\end{equation}
where we used (\ref{eq2.10c}) - (\ref{eq2.13}) and \cite[Eqs. (3.47) and (3.48)]{Dunster:2024:AZB}, and some calculation, to verify that the $\delta$ and $\delta^2$ terms cancel. Before continuing, we point out from (\ref{eq2.10b}) and (\ref{eq2.14}) that $F_{m}(\nu,\mathbf{z}_{\nu,m})=\mathcal{O}(\delta^{3})=\mathcal{O}(\nu^{-6})$ as $\delta \to 0$ ($\nu \to \infty$).

Next we construct lower and upper bounds in turn for the three terms in the braces on the RHS of (\ref{eq2.14}). Firstly, again by explicit differentiation and the chain rule, we have from (\ref{eq2.10a})
\begin{equation}
\label{eq2.15}
\breve{\zeta}'''(\delta''')
=\left(z_{m,1}+2\delta''' z_{m,2}\right)^{3}
\zeta'''(\mathbf{z}_{\nu,m})
+6 \left(z_{m,1}+2\delta''' z_{m,2}\right) z_{m,2} 
\zeta''(\mathbf{z}_{\nu,m}).
\end{equation}
But from \cref{lem:zeta} $\zeta''(z)$ and $-\zeta'''(z)$ are positive and strictly decreasing for $1 \leq z < \infty$. Furthermore, from \cref{lem:zmjMonotonic} $z_{m,1}+2\delta''' z_{m,2} \geq z_{m,1}+2 z_{m,2}>0$ for $1 \leq z_{m,0}<\infty$ and $\delta''' \in (0,1)$. Hence from (\ref{eq2.15}) and \cref{lem:bold(z[m])Bound}, and recalling $z_{m,2}<0$,
\begin{equation}
\begin{split}
\label{eq2.18}
& z_{m,1}^{3}
\zeta'''(z_{m,0})
+6  z_{m,1} z_{m,2}\zeta''(z_{m,0})
< \breve{\zeta}'''(\delta''')
\\
&< \left ( z_{m,1}+2z_{m,2}\right )^{3}
\zeta'''(z_{m,0}+\tfrac{1}{76})
+6 \left( z_{m,1}+2z_{m,2}\right) z_{m,2}
\zeta''(z_{m,0}+\tfrac{1}{76}) <0.
\end{split}
\end{equation}

Next, similarly to (\ref{eq2.10c}) and (\ref{eq2.15}),
\begin{equation*}
\breve{\Upsilon}_{1}''(\delta'')
=\left(z_{m,1}+2\delta'' z_{m,2}\right)^{2}
 \Upsilon_{1}''(\mathbf{z}_{\nu,m})
+2 z_{m,2} 
\Upsilon_{1}'(\mathbf{z}_{\nu,m}),
\end{equation*}
and so for $\delta'' \in (0,1)$ we have from \cref{lem:bold(z[m])Bound,lem:zmjMonotonic,lem:UpsilonMonotonic}
\begin{equation}
\begin{split}
\label{eq2.21}
0 & < \left(z_{m,1}+2 z_{m,2}\right)^{2}
\Upsilon_{1}''(z_{m,0}+\tfrac{1}{76})
+2 z_{m,2} \Upsilon_{1}'(z_{m,0}+\tfrac{1}{76})
\\ &
<  \breve{\Upsilon}_{1}''(\delta'')
< z_{m,1}^{2} \Upsilon_{1}''(z_{m,0})
+2 z_{m,2} \Upsilon_{1}'(z_{m,0}).
\end{split}
\end{equation}
Likewise
\begin{equation*}
\breve{\Upsilon}_{2}'(\delta)
=\left(z_{m,1}+2\delta z_{m,2}\right)
\Upsilon_{2}'(\mathbf{z}_{\nu,m}),
\end{equation*}
and thus for $\delta' \in (0,1)$
\begin{equation}
\label{eq2.23}
0 < \left(z_{m,1}+2z_{m,2}\right)
\Upsilon_{2}'(z_{m,0}+\tfrac{1}{76})
<\breve{\Upsilon}_{2}'(\delta')
< z_{m,1}\Upsilon_{2}'(z_{m,0}).
\end{equation}

For the last term on the RHS of (\ref{eq2.10b}), we have from \cref{lem:bold(z[m])Bound,lem:UpsilonMonotonic}
\begin{equation}
\label{eq2.24}
0 < \Upsilon_{3}(z_{m,0}+\tfrac{1}{76})
< \Upsilon_{3}(\mathbf{z}_{\nu,m})
< \Upsilon_{3}(z_{m,0}).
\end{equation}
Now from (\ref{eq2.10b}) and (\ref{eq2.14})
\begin{equation*}
F_{m}(\nu,\mathbf{z}_{\nu,m})
=\frac{1}{6\nu^{6}}
\left\{
\breve{\zeta}'''(\delta''')
+3\breve{\Upsilon}_{1}''(\delta'')
+6\breve{\Upsilon}_{2}'(\delta')
+6\Upsilon_{3}(\mathbf{z}_{\nu,m})
\right\},
\end{equation*}
and thus bearing in mind that $z_{m,1}=\hat{z}_{1}(z_{m,0})$ and $z_{m,2}=\hat{z}_{2}(z_{m,0})$ can be regarded as functions of $z_{m,0}$ we have from (\ref{eq2.18}) - (\ref{eq2.24})
\begin{equation}
\label{eq2.26}
\frac{\mathcal{F}_{1}(z_{m,0})}{6\nu^{6}}
\leq
F_{m}(\nu,\mathbf{z}_{\nu,m})
\leq 
\frac{\mathcal{F}_{2}(z_{m,0})}{6\nu^{6}},
\end{equation}
where
\begin{equation}
\begin{split}
\label{eq2.27}
\mathcal{F}_{1}(z_{m,0}) = &
z_{m,1}^{3}
\zeta'''(z_{m,0})
+6 z_{m,1} z_{m,2}\zeta''(z_{m,0})
\\ &
+3\left(z_{m,1}+2 z_{m,2}\right)^{2}
\Upsilon_{1}''(z_{m,0}+\tfrac{1}{76})
+6 z_{m,2} \Upsilon_{1}'(z_{m,0}+\tfrac{1}{76})
\\ &
+6\left(z_{m,1}+2z_{m,2}\right)
\Upsilon_{2}'(z_{m,0}+\tfrac{1}{76})
+6\Upsilon_{3}(z_{m,0}+\tfrac{1}{76}),
\end{split}
\end{equation}
and
\begin{equation}
\begin{split}
\label{eq2.28}
\mathcal{F}_{2}(z_{m,0})= &
\left ( z_{m,1}+2z_{m,2}\right )^{3}
\zeta'''(z_{m,0}+\tfrac{1}{76})
\\ &
+6 \left( z_{m,1}+2z_{m,2}\right) z_{m,2}
\zeta''(z_{m,0}+\tfrac{1}{76})
+3z_{m,1}^{2} \Upsilon_{1}''(z_{m,0})
 \\ &
+6 z_{m,2} \Upsilon_{1}'(z_{m,0})
+6 z_{m,1} \Upsilon_{2}'(z_{m,0})
+6\Upsilon_{3}(z_{m,0}).
\end{split}
\end{equation}

\begin{figure}[htbp]
 \centering
 \includegraphics[width=0.7\textwidth,
 keepaspectratio]
 {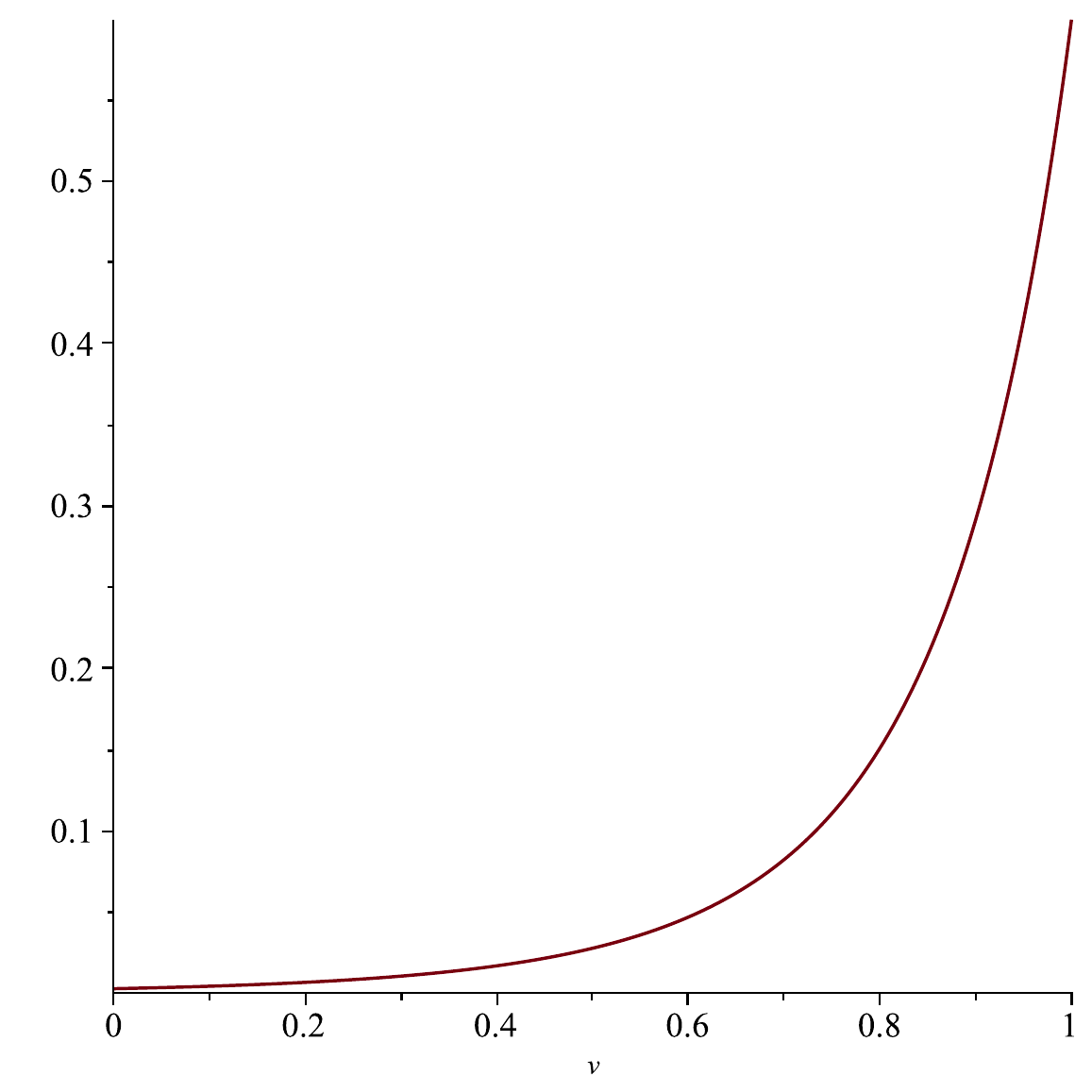}
 \caption{Graph of $p_{2}(v)$ for $0 \leq v < 1$}
 \label{fig:fig9}
\end{figure}

Consider next the denominator of the RHS of (\ref{eq2.5}). From \cref{lem:sigma,lem:UpsilonMonotonic} $\sigma(z)$ and $z|\eta'(1,z)|$ are decreasing for $1<z<\infty$, and hence so too is their product. Then we define
\begin{equation}
\label{eq2.7}
c_{1}:=1+\sup_{1 \leq z < \infty} \left\{
z \sigma(z) |\eta'(1,z)| \right\}
=1+\sigma(1)  |\eta'(1,1)|
= 1.0082524557 \cdots,
\end{equation}
where we utilised (\ref{eq1.9}) and (\ref{eq1.34}). Therefore from (\ref{eq1.35a}), (\ref{eq1.22}), (\ref{eq1.23}), (\ref{eq1.35}), (\ref{eq1.43}) and (\ref{eq2.7}) 
\begin{equation}
\label{eq2.6}
-\frac{c_{1}}{z \sigma(z)}
< F_{m}'(\nu,z) 
< -\frac{1}{z \sigma(z)} <0,
\end{equation}
for $1 \leq z < \infty$ and $1 \leq \nu < \infty$. 

Now from (\ref{eq1.3c}), (\ref{eq1.17}) - (\ref{eq1.19}),  (\ref{eq1.29}) - (\ref{eq1.31}) and (\ref{eq2.27})
\begin{equation}
\label{eq2.29}
\mathcal{F}_{1}(z_{m,0})
= \frac{6673\, (12^{2/3})}{58320\, \{z_{m,0}\}^{16/3}}
\left\{1+
\mathcal{O}\left(\frac{1}{z_{m,0}}\right)
\right\}
\quad (z_{m,0} \to \infty).
\end{equation}
In \cref{fig:fig9} we graph $p_{2}(v):=\{z_{m,0}\}^{16/3}\mathcal{F}_{1}(z_{m,0})$ with $z_{m,0}=1/(1-v)$ for $0 \leq v < 1$ (corresponding to $1 \leq z_{m,0} < \infty$). We introduced the factor $\{z_{m,0}\}^{16/3}$ to ensure that this function does not vanish at $v=1$ ($z=\infty$), which makes it clear that it is positive for all $z$.

Having confirmed that $\mathcal{F}_{1}(z_{m,0})$ is positive it is seen from (\ref{eq2.26}) that the same is true for $F_{m}(\nu,\mathbf{z}_{\nu,m})$. Consequently, from (\ref{eq2.5}) and (\ref{eq2.6}),
\begin{equation*}
0<\{c_{1}\}^{-1} \tau \sigma(\tau) 
F_{m}(\nu,\mathbf{z}_{\nu,m})
<\mathbf{e}_{\nu,m} 
<  \tau \sigma(\tau) 
F_{m}(\nu,\mathbf{z}_{\nu,m}),
\end{equation*}
and so from \cref{lem:sigma} and (\ref{eq2.4}),
\begin{multline}
\label{eq2.9}
\{c_{1}\}^{-1} z_{m,0}\sigma(z_{m,0})
F_{m}(\nu,\mathbf{z}_{\nu,m})
<  \mathbf{e}_{\nu,m} 
\\
< \left(z_{m,0}+\tfrac{1}{73}\right) 
\sigma\left(z_{m,0}+\tfrac{1}{73}\right) 
F_{m}(\nu,\mathbf{z}_{\nu,m}).
\end{multline}

\begin{figure}[htbp]
 \centering
 \includegraphics[width=0.7\textwidth,
 keepaspectratio]
 {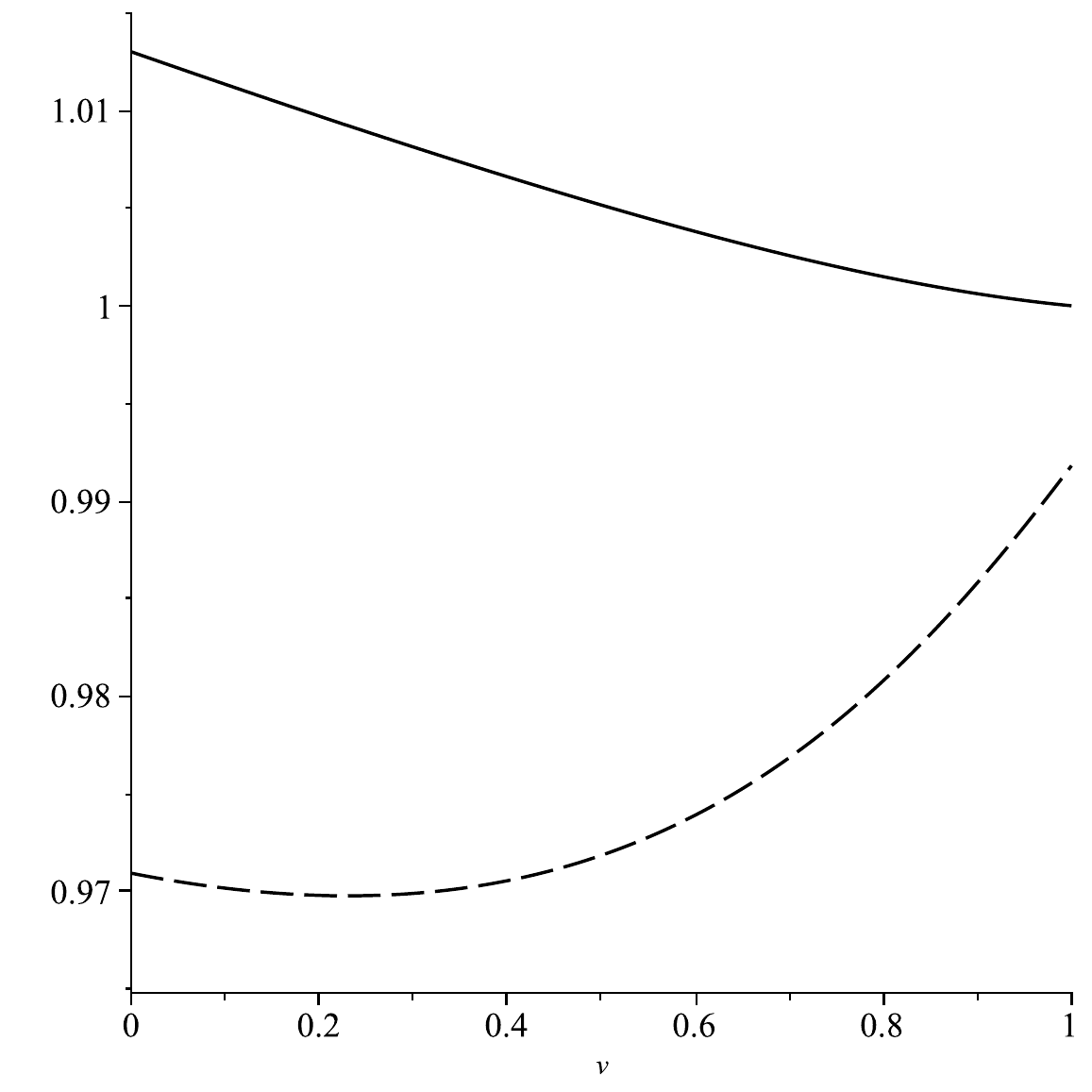}
 \caption{Graphs of $\mathcal{G}_{1}((1-v)^{-1})$ (dashed) and $\mathcal{G}_{2}((1-v)^{-1})$ (solid) for $0 \leq v < 1$}
 \label{fig:fig10}
\end{figure}

At this stage we could insert the bounds (\ref{eq2.26}) to obtain explict lower and upper bounds for $\mathbf{e}_{\nu,m}$, but we prefer to go one step further and simplify these bounds to only involve $\nu$ and $z_{m,3}$, at a small sacrifice in sharpness. To this end, we have from (\ref{eq2.26}) and (\ref{eq2.9})
\begin{equation}
\label{eq2.30}
\frac{\mathcal{G}_{1}(z_{m,0}) z_{m,3}}{\nu^6}
<\mathbf{e}_{\nu,m} 
< \frac{\mathcal{G}_{2}(z_{m,0}) z_{m,3}}{\nu^6},
\end{equation}
where
\begin{equation}
\label{eq2.31}
\mathcal{G}_{1}(z_{m,0})=
\frac{z_{m,0}\sigma(z_{m,0})
\mathcal{F}_{1}(z_{m,0})}{6 c_{1} z_{m,3}},
\end{equation}
and
\begin{equation}
\label{eq2.32}
\mathcal{G}_{2}(z_{m,0})=
\frac{(z_{m,0}+\tfrac{1}{73})\sigma(z_{m,0}
+\tfrac{1}{73})\mathcal{F}_{2}(z_{m,0})}{6 z_{m,3}},
\end{equation}
again remembering that $z_{m,3}=\hat{z}_{3}(z_{m,0})$ can be regarded as a function of $z_{m,0}$.

The graphs of $\mathcal{G}_{1}((1-v)^{-1})$ (dashed curve) and $\mathcal{G}_{2}((1-v)^{-1})$ (solid curve) for $0\leq v <1$ (equivalent to $\mathcal{G}_{j}(z_{m,0})$ ($j=1,2$) for $1 \leq z_{m,0} < \infty$) are depicted in \cref{fig:fig10}. From this $\mathcal{G}_{1}((1-v)^{-1})$ has a minimum value found to be $0.9697464085\cdots$ at $v = 0.2307692292\cdots$ (approximately $z_{m,0}=1.3$). 

Further, $\mathcal{G}_{2}((1-v)^{-1})$ is seen to be monotonically decreasing, and from (\ref{eq1.3b}), (\ref{eq1.9}), (\ref{eq1.14}) - (\ref{eq1.16}), (\ref{eq1.27}) - (\ref{eq1.28a}), (\ref{eq2.28}), (\ref{eq2.7}) and (\ref{eq2.32}) the maximum value is given by $\mathcal{G}_{2}(1)= 1.0130228266 \cdots$. We also note in passing that from (\ref{eq1.3c}), (\ref{eq1.8}) and (\ref{eq1.17}) - (\ref{eq1.19}) $\lim_{z_{m,0} \to \infty}\mathcal{G}_{2}(z_{m,0})= 1$.

From (\ref{eq2.30}) and these extrema of $\mathcal{G}_{j}(z_{m,0})$ ($j=1,2$) we have established (\ref{eq1.48}), and the proof of the theorem is complete.

\appendix
\section{Proofs of \cref{lem:zeta,lem:sigma,lem:zmjMonotonic,lem:bold(z[m])Bound,lem:UpsilonMonotonic,lem:SumUpsilon,lem:ScriptZ,lem:etaBound,lem:zhatbound,lem:wnumBounds,lem:MBound}}
\label{secA}

\begin{proof}[Proof of \cref{lem:zeta}]
That all the derivatives of $\zeta$ approach zero as $z \to \infty$ follows from (\ref{eq1.3c}). From (\ref{eq1.35a}) and repeated use of the chain rule the fourth derivative of $\zeta$ with respect to $z$ is given explicitly by
\begin{equation}
\label{eqA.32aa}
\zeta^{(4)}
= \frac {6}{z^{4} \sigma}
+\frac {6 z^{2} \sigma^{3}-11}{2
z^{4} \sigma^{2} \zeta}
+\frac {6-5 z^{2} \sigma^{3}}
{z^{4} \sigma^{3} \zeta^{2}}
+\frac {\left( 2 z^{2} \sigma^{3}-1\right)  
\left( 3 z^{4} \sigma^{6} 
+2 z^{2} \sigma^{3}+7\right) }{2 z^{4} \sigma^{4} \zeta^{3}}.
\end{equation}

\begin{figure}[htbp]
 \centering
 \includegraphics[width=0.7\textwidth,
 keepaspectratio]
 {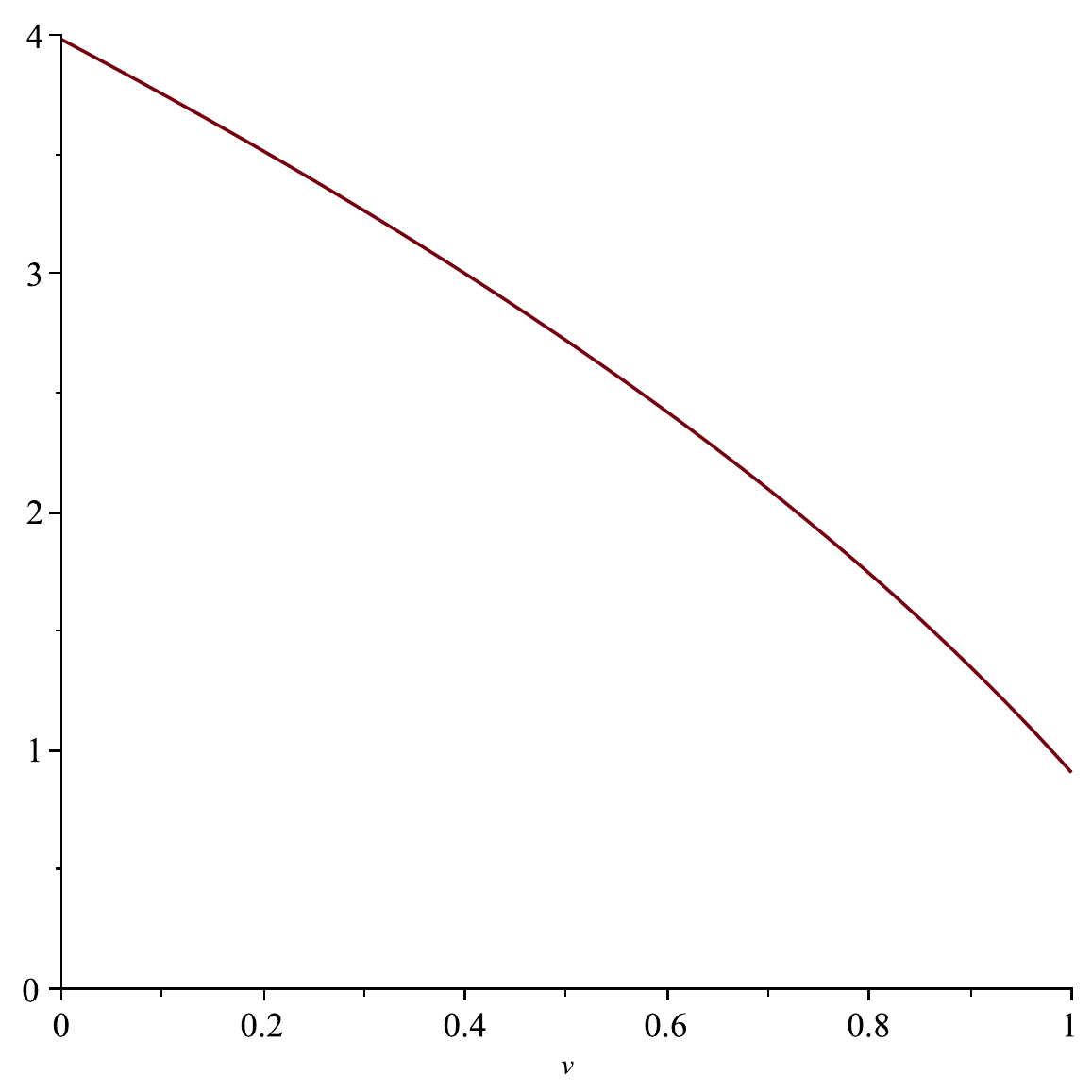}
 \caption{Graph of $(1-v)^{-10/3}\zeta^{(4)}((1-v)^{-1})$ for $0 \leq v < 1$}
 \label{fig:figZeta}
\end{figure}

In \cref{fig:figZeta} we graph $z^{10/3}\zeta^{(4)}(z)$, with $z$ replaced by $1/(1-v)$ for $0 \leq v < 1$, which corresponds to $1 \leq z < \infty$. The factor $z^{10/3}=(1-v)^{-10/3}$ ensures that the function does not vanish at $v=1$ ($z=\infty$), noting from (\ref{eq1.3c}) that $\zeta^{(4)}(z)=\mathcal{O}(z^{-10/3})$ as $z \to \infty$. From this we deduce that $\zeta^{(4)}(z)$ is positive, and recall from (\ref{eq1.3c}) that it approaches zero (through positive values) as $v \to 1^{-}$ ($z \to \infty$).

By integration with respect to $z$ it follows that $\zeta'''$ is monotonically increasing, and must be negative from its (negative) value at $z=0$ (see (\ref{eq1.3b})), its limit being zero through negative values as $z \to \infty$ (again from (\ref{eq1.3c})), and by continuity. Thus $-\zeta'''$ is positive and decreases monotonically to zero, as asserted. Integrating twice more with similar arguments yields the stated results for $\zeta''$ and $-\zeta'$.
\end{proof}

\begin{proof}[Proof of \cref{lem:zmjMonotonic}]
The plot of the functions
\begin{equation}
\label{eqA.33z}
\begin{split}
p_{3}(v)&:=-\lambda_{3}(v)(1-v)^{-2}
\hat{z}_{1}'((1-v)^{-1})
\quad  (\text{solid curve}),
\\
p_{4}(v)&:=\lambda_{4}(v)(1-v)^{-4}
\hat{z}_{2}'((1-v)^{-1})
\quad  (\text{dashed curve}),
\\
p_{5}(v)&:=-\lambda_{5}(v)(1-v)^{-6}
\hat{z}_{3}'((1-v)^{-1})
\quad  (\text{dotted curve}),
\\
p_{6}(v)&:=\lambda_{6}(v)(1-v)^{-1}
\left\{\hat{z}_{1}((1-v)^{-1})
+2\hat{z}_{2}((1-v)^{-1})\right\}
\quad  (\text{dash-dotted curve}),
\end{split}
\end{equation}
is shown in \cref{fig:fig1} for $0 \leq v <1$ (which is equivalent to $1 \leq z < \infty$). For $j=1,2,3$ the factor $(1-v)^{-2j}$ was introduced so that $p_{j}(v)$ ($j=3,4,5$) do not vanish at $v=1$ ($z=\infty$), since from (\ref{eq1.17}) - (\ref{eq1.19}) we see that $\hat{z}_{j}'(z)=\mathcal{O}(z^{-2j})$ as $z \to \infty$. Similarly for the factor $(1-v)^{-1}$ in $p_{6}(v)$. Each $\lambda_{j}(v)$ is a positive scaling factor introduced for convenience so that $p_{j}(0)=\lim_{v \to 1^{-}}p_{j}(v)=1$, and we chose them to be of the form 
\begin{equation}
\label{eqA.33l}
\lambda_{j}(v)
=\alpha_{j}\beta_{j}^{v}.
\end{equation}
So, for example, we have from (\ref{eq1.14}) and  (\ref{eq1.17})
\begin{equation*}
\alpha_{3}=
-\{\hat{z}_{1}'(1)\}^{-1}
=3150, \quad
\beta_{3}=-\left[\alpha_{3}
\lim_{z \to \infty}\left\{z^{2}
\hat{z}_{1}'(z)\right\}\right]^{-1}
=\tfrac{1}{175}.
\end{equation*}
Similar scaling factors are used in various functions below, namely those defined by (\ref{eqA.06b}) and (\ref{eqA.06c}).

\begin{figure}[htbp]
 \centering
 \includegraphics[width=0.7\textwidth,
 keepaspectratio]
 {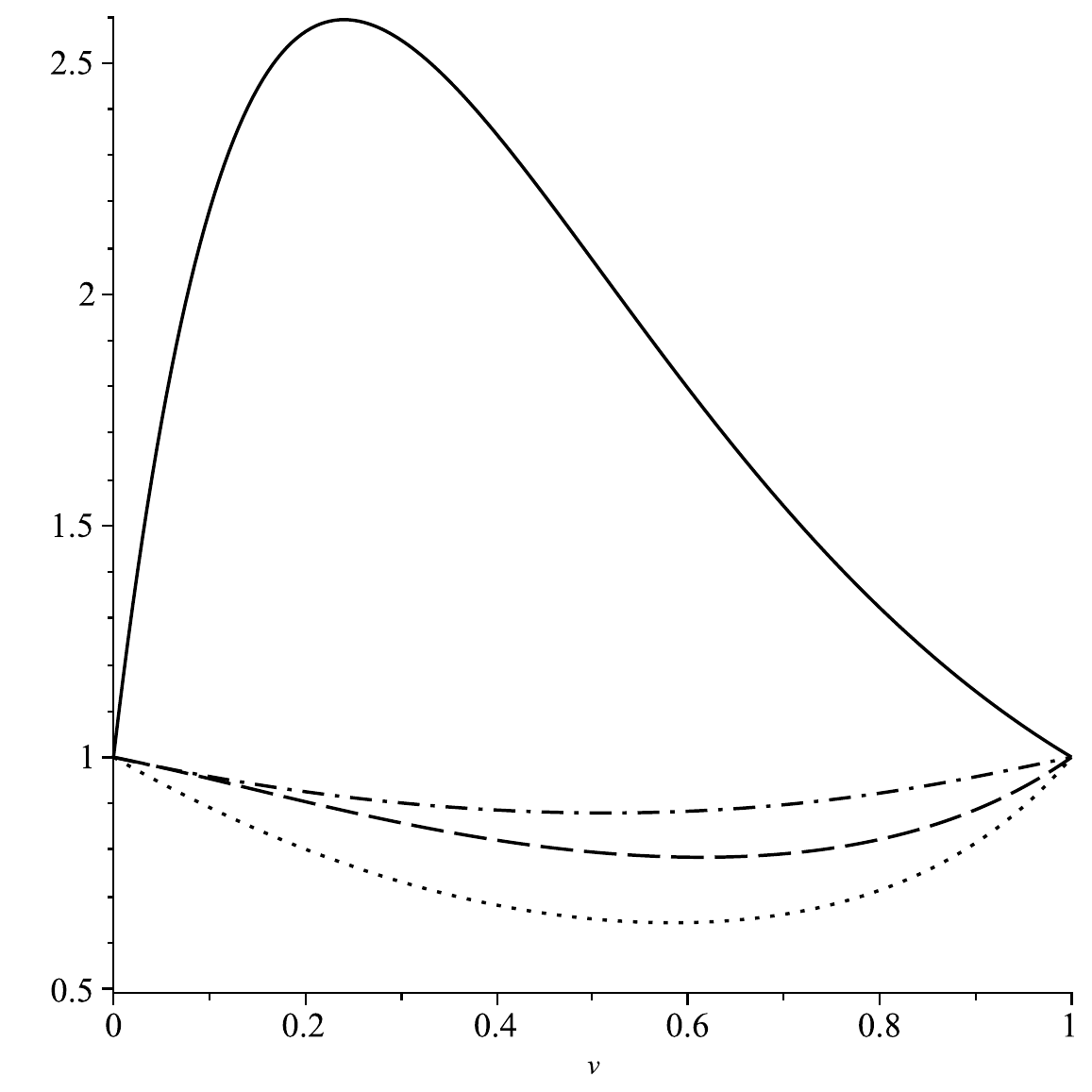}
 \caption{Graphs of $p_{3}(v)$ (solid), $p_{4}(v)$ (dashed), $p_{5}(v)$ (dotted), and $p_{6}(v)$ (dash-dotted) for $0 \leq v <1$}
 \label{fig:fig1}
\end{figure}

We observe from their graphs that $p_{j}(v)>0$ ($j=3,4,5,6$) for $0 \leq v <1$, the first three of which imply $(-1)^{j-1}\hat{z}_{j}'(z) <0$ ($j=1,2,3$) for $1 \leq z < \infty$. Consequently, in conjunction with (\ref{eq1.14}) - (\ref{eq1.19}), it follows on integration that each $(-1)^{j-1}\hat{z}_{j}(z)$ is positive and monotonically decreasing to 0 for $1 \leq z < \infty$. 

Finally, $p_{6}(v)>0$ ($0 \leq v <1$) and $\hat{z}_{2}(z)<0$ ($1 \leq z<\infty$) establishes (\ref{eq1.20}).
\end{proof}

\begin{proof}[Proof of \cref{lem:bold(z[m])Bound}]
From \cref{lem:zmjMonotonic} we have $z_{m,2}<0$ and $z_{m,1}+2z_{m,2}>0$ for $1 <  z_{m,0} < \infty$, recalling that $z_{m,j}=\hat{z}_{j}(z_{m,0})$ ($j=1,2$). Now let $\delta = \nu^{-2}$ and define $h_{m}(\delta):=\delta z_{m,1}+\delta^{2} z_{m,2}$ for $0 < \delta  \leq 1$. Note from (\ref{eq1.12}) that $\mathbf{z}_{\nu,m}=z_{m,0}+h_{m}(\delta)$. Now $h_{m}'(\delta)=z_{m,1}+2\delta z_{m,2} \geq z_{m,1}+2z_{m,2}>0$, and consequently $0=h_{m}(0) < h_{m}(\delta) \leq h_{m}(1)=z_{m,1}+z_{m,2}$. This gives
\begin{equation}
\label{eqA.03}
z_{m,0} < \mathbf{z}_{\nu,m} 
\leq z_{m,0}
+z_{m,1} + z_{m,2},
\end{equation}
and as such the lower bound of (\ref{eq1.21}) is proven.

\begin{figure}[htpb]
 \centering
 \includegraphics[width=0.7\textwidth,
 keepaspectratio]
 {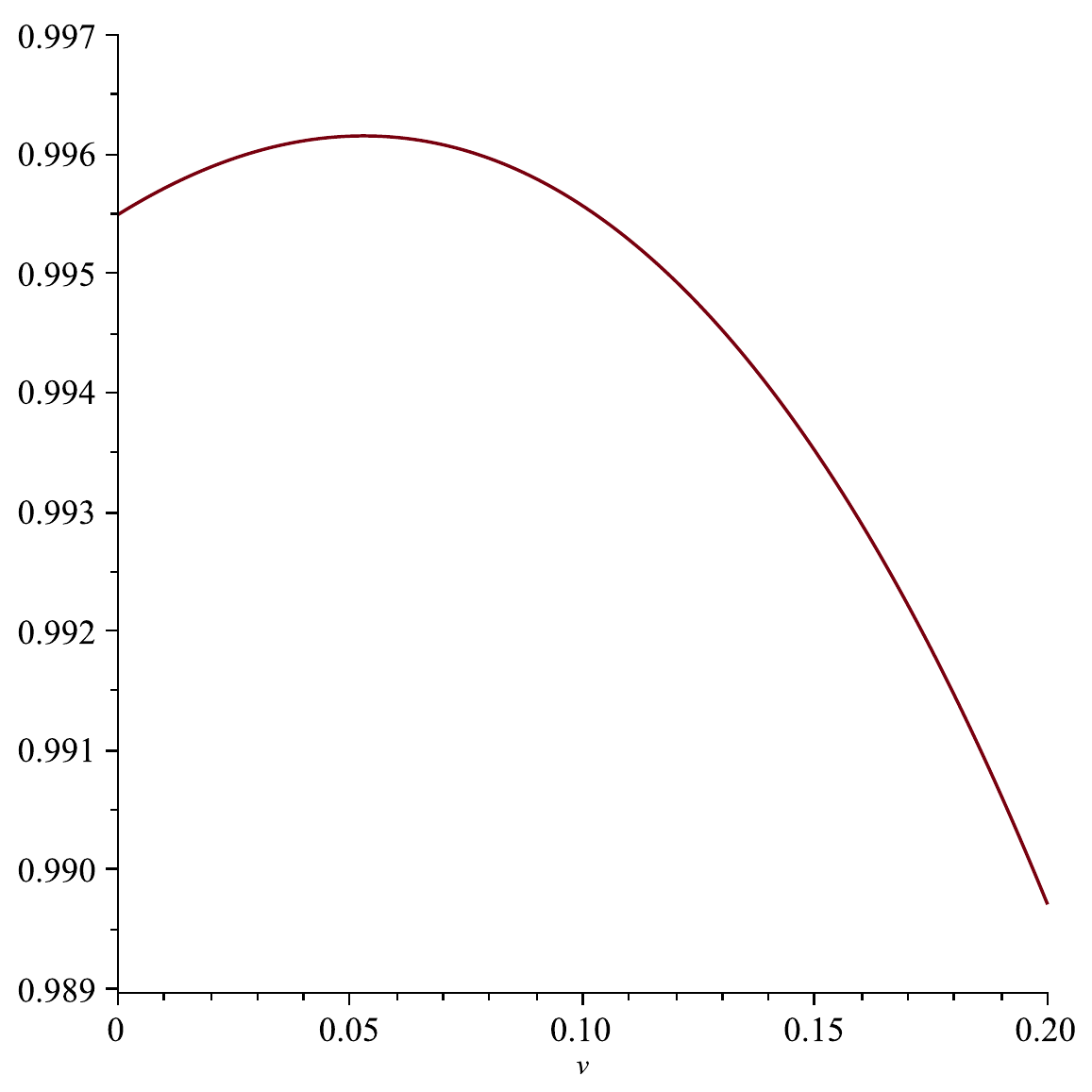}
 \caption{Graph of $p_{7}(v)$ for $0 \leq v \leq 0.2$}
 \label{fig:fig13}
\end{figure}

Next we find numerically that $p_{7}(v):=76\{\hat{z}_{1}(z)+\hat{z}_{2}(z)\}$ with $z=1/(1-v)$ attains for $v \in [0,1)$ a maximum value of $0.99615\cdots$ at $v = 0.05288\cdots$. A graph of $p_{7}(v)$ for $0 \leq v \leq 0.2$ is depicted in \cref{fig:fig13}, and we also find it decreases to zero for $0.2 \leq v <1$. We conclude that $76(z_{m,1}+z_{m,2})<1$  for $1 \leq z_{m,0} < \infty$, and accordingly it follows from (\ref{eqA.03}) that the upper bound of (\ref{eq1.21}) is true.
\end{proof}

\begin{proof}[Proof of \cref{lem:UpsilonMonotonic}]
Firstly, $\Upsilon_{1}(z)$, $-\Upsilon_{2}(z)$, $\Upsilon_{3}(z)$, $-\Upsilon_{1}'(z)$, $\Upsilon_{2}'(z)$, $-\Upsilon_{3}'(z)$,  $\Upsilon_{1}''(z)$, $\eta(1,z)$, $-\eta'(1,z)$, $-z\eta'(1,z)$ and $z\sigma(z)\eta(1,z)$ all approach zero through positive values as $z \to \infty$, as can be confirmed from (\ref{eq1.9b}), (\ref{eq1.23}) and (\ref{eq1.29}) - (\ref{eq1.31}).

\begin{figure}[htbp]
 \centering
 \includegraphics[width=0.7\textwidth,
 keepaspectratio]
 {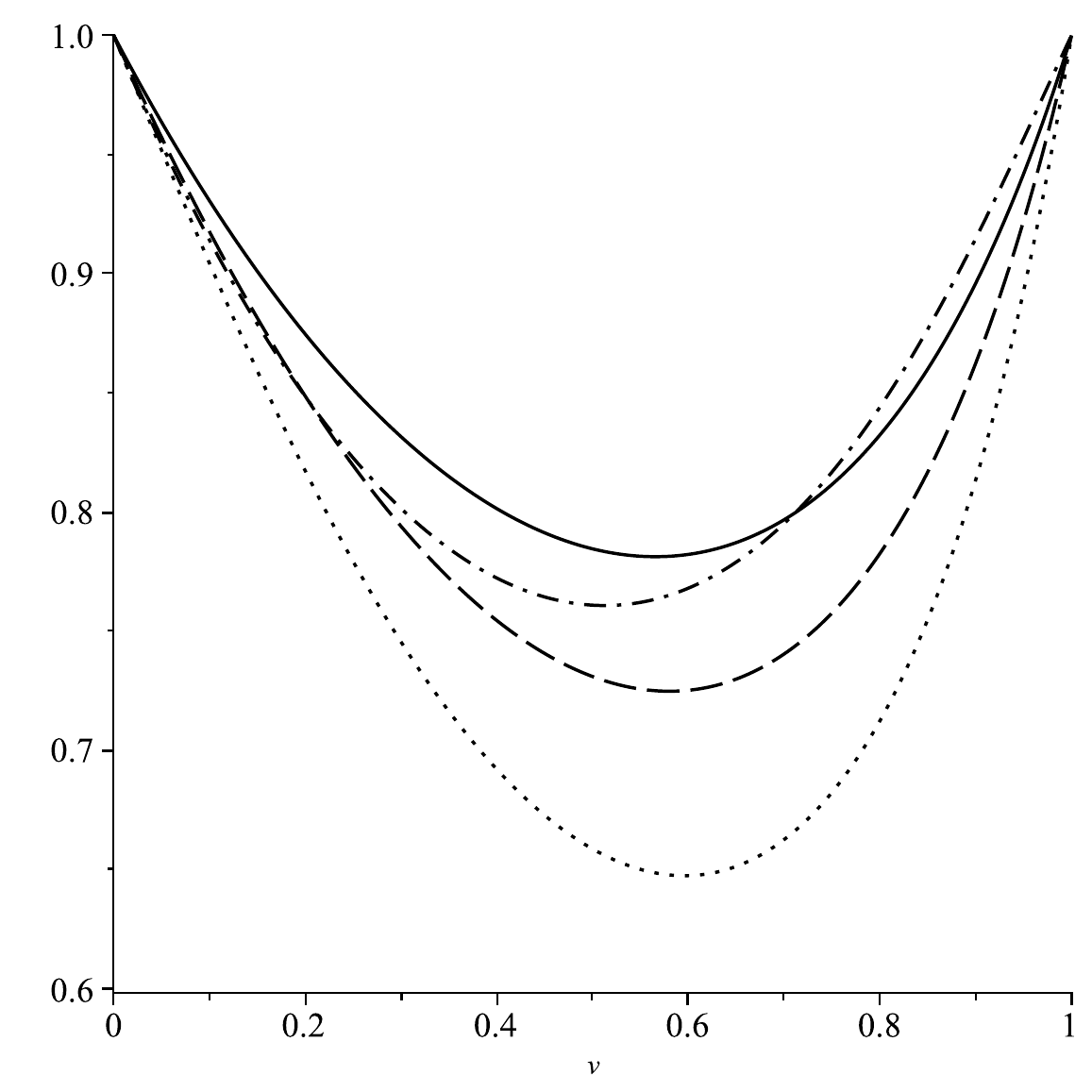}
 \caption{Graphs of $p_{8}(v)$ (solid), $p_{9}(v)$ (dashed), $p_{10}(v)$ (dotted), and $p_{11}(v)$ (dash-dotted) for $0 \leq v <1$}
 \label{fig:fig2}
\end{figure}

\begin{figure}[htbp]
 \centering
 \includegraphics[width=0.7\textwidth,
 keepaspectratio]
 {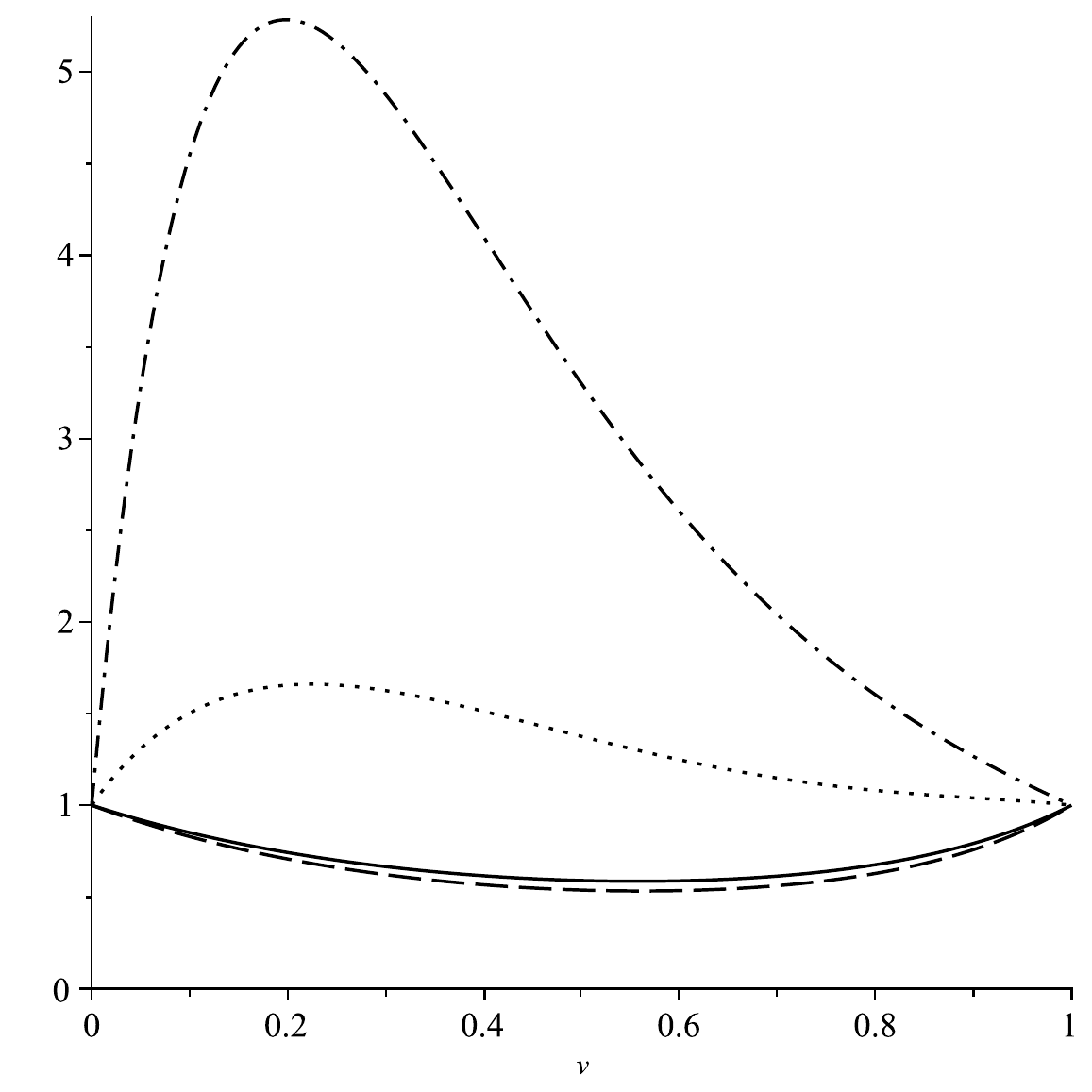}
 \caption{Graphs of $p_{12}(v)$ (solid), $p_{13}(v)$ (dashed), $p_{14}(v)$ (dotted), and $p_{15}(v)$ (dash-dotted) for $0 \leq v <1$}
 \label{fig:fig2a}
\end{figure}

Next, in \cref{fig:fig2} for $0 \leq v <1$ we plot the following functions which are shown to all be positive:
\begin{equation}
\label{eqA.06b}
\begin{split}
p_{8}(v)&:=-\lambda_{8}(v)
(1-v)^{-13/3}\Upsilon_{1}'''((1-v)^{-1})
\quad  (\text{solid curve}),
\\
p_{9}(v)&:=-\lambda_{9}(v)(1-v)^{-16/3}
\Upsilon_{2}''((1-v)^{-1})
\quad  (\text{dashed curve}),
\\
p_{10}(v)&:=\lambda_{10}(v)(1-v)^{-22/3}
\Upsilon_{3}''((1-v)^{-1})
\quad  (\text{dotted curve}),
\\
p_{11}(v)&:=\lambda_{11}(v) (1-v)^{-10/3}
\textstyle \sum_{j=1}^{2}\Upsilon_{j}''((1-v)^{-1})
\quad  (\text{dash-dotted curve}).
\end{split}
\end{equation}
Similarly to (\ref{eqA.33z}), $\lambda_{j}(v)$ are again positive scaling factors of the form (\ref{eqA.33l}), with the two constants therein being chosen so that all the functions take the value 1 at the end points.

From \cref{fig:fig2a} for $0 \leq v <1$ the following functions are also all seen to be positive:
\begin{equation}
\label{eqA.06c}
\begin{split}
p_{12}(v)&:=-\lambda_{12}(v)(1-v)^{-20/3}
\left\{\Upsilon_{2}^{2}((1-v)^{-1})
\right.
\\ &
\left.
-\tfrac{64}{25}\Upsilon_{1}((1-v)^{-1})
\Upsilon_{3}((1-v)^{-1})\right\}
\quad   (\text{solid curve}),
\\
p_{13}(v)&:=-\lambda_{13}(v)(1-v)^{-20/3}
\left\{\Upsilon_{2}'^{2}((1-v)^{-1})
\right.
\\ &
\left.
-3\Upsilon_{1}'((1-v)^{-1})
\Upsilon_{3}'((1-v)^{-1})\right\}
\quad  (\text{dashed curve}),
\\ p_{14}(v)&:=\lambda_{14}(v)
\left[z^{7/3}
\left\{z\eta'(1,z)\right\}'\right]
_{z=(1-v)^{-1}}
\quad  (\text{dotted curve}),
\\ p_{15}(v)&:= -\lambda_{15}(v)  \left[
z^{2}\left\{z \sigma(z) \eta(1,z)\right\}'\right]
_{z=(1-v)^{-1}}
\quad  (\text{dash-dotted curve}).
\end{split}
\end{equation}

Note since $p_{10}(v)>0$ and $p_{11}(v)>0$ that
\begin{equation}
\label{eqA.10c}
\eta''(1,z)
>\Upsilon_{1}''(z)+\Upsilon_{2}''(z)>0
\quad  (1 \leq z < \infty).
\end{equation}
The stated results then follow by integration of all of the above, along with their behaviour at $z=\infty$, similarly to the proof of \cref{lem:zmjMonotonic}.
\end{proof}

\begin{proof}[Proof of \cref{lem:SumUpsilon}]
From \cref{lem:UpsilonMonotonic} we have $\eta(1,z)$ is (strictly) decreasing, this gives the third inequality of (\ref{eq1.32}). Again let $\delta=\nu^{-2}$ and consider 
\begin{equation}
\label{eqA.10}
h(\delta,z):=\eta(\delta^{-1/2},z)
=\Upsilon_{1}(z)\delta
+\Upsilon_{2}(z)\delta^{2}+\Upsilon_{3}(z)\delta^{3}.
\end{equation}
Then $\partial h/\partial \delta =0$ for 
\begin{equation*}
\delta= \frac{|\Upsilon_{2}(z)|
\pm \left\{\Upsilon_{2}^{2}(z)
-3\Upsilon_{1}(z)\Upsilon_{3}(z)\right\}^{1/2}}
{3\Upsilon_{3}(z)}.
\end{equation*}
Now from (\ref{eqA.18}) and the positivity of $\Upsilon_{1}(z)\Upsilon_{3}(z)$
\begin{equation}
\label{eqA.13}
\Upsilon_{2}^{2}(z)-3\Upsilon_{1}(z)\Upsilon_{3}(z)<\Upsilon_{2}^{2}(z)
-\tfrac{64}{25}\Upsilon_{1}(z)\Upsilon_{3}(z)<0.
\end{equation}
Thus  $\partial h/\partial \delta \neq 0$ for $1 \leq z < \infty$ and $0 \leq \delta \leq 1$. Moreover $h(0,z)=0$ and $h(1,z)=\eta(1,z)>0$ ($1 < z < \infty$), and this latter value then is the maximum value of $h(\delta,z)$ for $0 \leq \delta \leq 1$ ($1 \leq \nu < \infty$) for each $z \in (1,\infty)$. This establishes the second inequality of (\ref{eq1.32}).

Next from (\ref{eqA.10c}) and (\ref{eqA.10}) for $0 < \delta \leq 1$, recalling $\Upsilon_{2}''(z)<0$ and $\Upsilon_{3}''(z)>0$,
\begin{equation*}
\delta^{-1}h''(\delta,z)=\Upsilon_{1}''(z)
+\Upsilon_{2}''(z)\delta+\Upsilon_{3}''(z)\delta^{2}
>\Upsilon_{1}''(z)+\Upsilon_{2}''(z)>0.
\end{equation*}
Therefore from (\ref{eq1.23}), (\ref{eq1.29}) - (\ref{eq1.31}) and (\ref{eqA.10}) it follows that $h(\delta,z)$  and $-h'(\delta,z)$ are positive and strictly decreasing to zero for $1 \leq z < \infty$. The same of course is true for $\eta(\nu,z)$  and $-\eta'(\nu,z)$. This gives the first (strict) inequality of (\ref{eq1.32}), as well as the last (strict) inequality of (\ref{eq1.34}).

For (\ref{eq1.33}) consider
\begin{equation*}
\hat{h}(\delta,z):=\delta^{-1/3}h(\delta,z)
=\Upsilon_{1}(z)\delta^{2/3}
+\Upsilon_{2}(z)\delta^{5/3}+\Upsilon_{3}(z)\delta^{8/3},
\end{equation*}
and as such
\begin{equation*}
\frac{\partial \hat{h}(\delta,z)}
{\partial \delta}=
\frac{1}{3\delta^{1/3}}\left\{
2\Upsilon_{1}(z) + 5\Upsilon_{2}(z)\delta
+8\Upsilon_{3}(z)\delta^{2}\right\}.
\end{equation*}
Note this derivative approaches $+\infty$ as $\delta \to 0^{+}$ on account of (\ref{eq1.27}). From (\ref{eqA.18}) we similarly find that as a function of $\delta \in (0,1]$ $\hat{h}(\delta,z)$ has no real zeros for each $z \in [1, \infty)$, and so is positive for all such values of these variables. We deduce in a similar manner to above for each fixed $z$ that $\hat{h}(\delta,z)$ is positive and increasing for $0 < \delta \leq 1$, and (\ref{eq1.33}) follows.

Next the leading equality of (\ref{eq1.34}) comes from (\ref{eq1.27}) - (\ref{eq1.28a}). In addition, from \cref{lem:UpsilonMonotonic} $\eta'(1,z)$ is (strictly) increasing, and consequently the first inequality of (\ref{eq1.34}) follows.

Finally, again let $\delta=\nu^{-2}$, and consider 
\begin{equation}
\label{eqA.20}
h'(\delta,z)=\partial h(\delta,z) /dz
=\Upsilon_{1}'(z)\delta
+\Upsilon_{2}'(z)\delta^{2}+\Upsilon_{3}'(z)\delta^{3},
\end{equation}
and then $\partial h'/\partial \delta =0$ for 
\begin{equation*}
\delta= \frac{\Upsilon_{2}'(z)
\pm \left\{\Upsilon_{2}'^{2}(z)
-3\Upsilon_{1}'(z)\Upsilon_{3}'(z)\right\}^{1/2}}
{3|\Upsilon_{3}'(z)|}.
\end{equation*}
Therefore from (\ref{eqA.22}) we deduce that $\partial h'/\partial \delta \neq 0$ for $1 < z < \infty$ and $0 \leq \delta \leq 1$. However from (\ref{eqA.20}) $h'(0,z)=0$, and $h'(1,z)<0$ as we have just shown, and this latter value then is the minimum value of $h'(\delta,z)$ for each $0 \leq \delta \leq 1$ ($1 \leq \nu < \infty$) and for each $z \in (1,\infty)$. This establishes the second inequality of (\ref{eq1.34}), and the proof of the lemma is complete.
\end{proof}

\begin{proof}[Proof of \cref{lem:ScriptZ}]

For $1 \leq z < \infty$ and $1 \leq \nu < \infty$ we have from (\ref{eq1.35a}), (\ref{eq1.22}), (\ref{eq1.23}) and (\ref{eq1.34})
\begin{equation*}
\mathcal{Z}_{3}'(\nu,z)
=-\frac{1}{z \sigma(z)}
+ \eta'(\nu,z)
< -\frac{1}{z \sigma(z)},
\end{equation*}
and (\ref{eq1.35}) then follows, noting that $\sigma(z)>0$ for $1 \leq z<\infty$.

Finally, from (\ref{eq1.5}), (\ref{eq1.22}), (\ref{eq1.23}),  (\ref{eq1.27}) - (\ref{eq1.28a}), (\ref{eq1.33}), (\ref{eq1.35}), \cref{lem:UpsilonMonotonic}, and noting $\mathrm{a}_{m} \leq \mathrm{a}_{1}=-2.338107410 \cdots$, for $z \geq z_{m,0}>1$ ($\zeta \leq \zeta_{m,0} <0$) we have
\begin{equation*}
\begin{split}
\nu^{2/3} \mathcal{Z}_{3}(\nu,z) & \leq 
\nu^{2/3} \zeta_{m,0}
+\nu^{2/3} \eta(\nu,z_{m,0})
< \mathrm{a}_{1} + \nu^{2/3} \eta(\nu,1)
\\
&  < -2.338
+ \nu^{2/3} \eta(\nu,1)
\leq -2.338 + \eta(1,1)
\\
&
=-2.338 +\tfrac{44873962351}
{3302530481250}2^{1/3} <0.
\end{split}
\end{equation*}
Hence $\mathcal{Z}_{3}(\nu,z)<0$ for $z_{m,0} \leq z < \infty$ ($-\infty<\zeta \leq \zeta_{m,0}$) and $\nu \geq 1$. 
\end{proof}

\begin{proof}[Proof of \cref{lem:wnumBounds}]

Since $\zeta =0$ when $z=1$ we see from (\ref{eq1.22}), (\ref{eq1.23}) and (\ref{eq1.33}) that $\mathcal{Z}_{3}(\nu,1)>0$. On the other hand from (\ref{eq1.3c}), (\ref{eq1.22}), (\ref{eq1.23}), (\ref{eq1.29}) - (\ref{eq1.31}) it is seen that $\mathcal{Z}_{3}(\nu,z) \to -\infty$ as $z \to \infty$ for $1 \leq \nu < \infty$. Thus from \cref{lem:ScriptZ} $\mathcal{Z}_{3}(\nu,z)$ decreases monotonically from a positive value to $-\infty$ for $1 \leq z < \infty$. But from (\ref{eq1.51}) and (\ref{eq1.57})
\begin{equation*}
\nu^{2/3} \mathcal{Z}_{3}(\nu,w_{\nu,1}^{-})
=\mathrm{a}_{1}+r_{1}^{+}
=-\tfrac{299}{800}\left(24\pi^{2}\right)^{1/3}<0,
\end{equation*}
confirming that $w_{\nu,1}^{-}>1$.

Next define (see (\ref{eq1.50}), (\ref{eq1.52}) and (\ref{eq1.57}))
\begin{equation}
\label{eqA.41}
\alpha_{\nu,m}^{-}
:=\mathcal{Z}_{3}(\nu,w_{\nu,m}^{-})
=\frac{\mathrm{a}_{m,0}}{\nu^{2/3}}
\left(1-\frac{0.01}{4m-1}\right)
\quad  (m=1,2,3,\ldots).
\end{equation}
We note that $|\mathrm{a}_{m}|>|\mathrm{a}_{m,0}|$ for all $m$ (see \cite{Pittaluga:1991:IZA}), and hence from (\ref{eq1.5}) and (\ref{eqA.41})
\begin{equation}
\label{eqA.42}
\frac{\alpha_{\nu,m}^{-}}{\zeta_{m,0}}
=\frac{\mathrm{a}_{m,0}}{\mathrm{a}_{m}}
\left(1-\frac{0.01}{4m-1}\right)<1
\quad  (m=1,2,3,\ldots).
\end{equation}

A lower bound is the goal, and from explicit computation of the middle term of (\ref{eqA.42}) for $m=1$ we find that $\alpha_{\nu,1}^{-}/\zeta_{1,0}=0.98905\cdots > 0.989$. Now, from the equality in (\ref{eqA.42}) and \cite[Eqs. (4.4) and (4.5)]{Qu:1999:BPU} (which are proven in \cite{Hethcote:1970:EBA} and \cite{Wong:1991:OPI}), a general bound is given by
\begin{equation}
\label{eqA.43}
\frac{\alpha_{\nu,m}^{-}}{\zeta_{m,0}}
\geq
\left\{1+\frac{0.130}{\left(
\frac{3}{8}\pi(4m-1.051)\right)^{2}}\right\}^{-1}
\left(1-\frac{0.01}{4m-1}\right)
\quad  (m=1,2,3,\ldots).
\end{equation}
For $m=2,3,4,\ldots$ the RHS is found using elementary calculus to attain an absolute minimum at $m=2$ of $0.99663\cdots > 0.989$. Thus we have for $m=1,2,3,\ldots$, $1 \leq \nu < \infty$
\begin{equation}
\label{eqA.44}
0.989|\zeta_{m,0}|<|\alpha_{\nu,m}^{-}|<|\zeta_{m,0}|.
\end{equation}

\begin{figure}[htpb]
 \centering
 \includegraphics[width=0.7\textwidth,
 keepaspectratio]
 {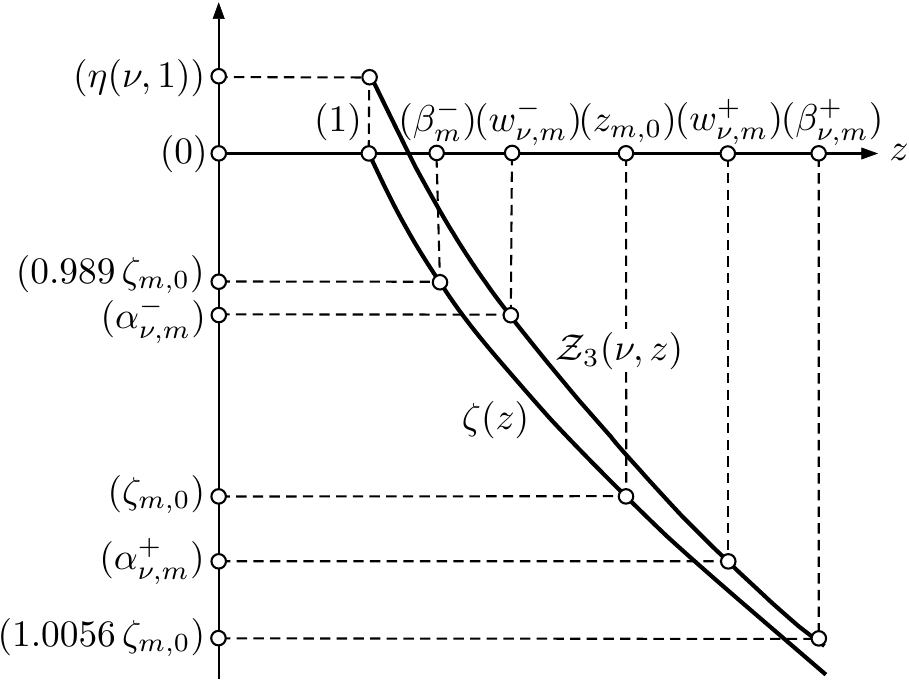}
 \caption{Graph illustrating $\beta_{m}^{-}<
 w_{\nu,m}^{-}<w_{\nu,m}^{+}<\beta_{\nu,m}^{+}$}
 \label{fig:zeta-Z}
\end{figure}

Next let $\beta_{m}^{-} =\zeta^{-1}(0.989\,\zeta_{m,0})$ (here and below we mean the inverse of $\zeta$). Then from (\ref{eqA.44}) and \cref{fig:zeta-Z} (not drawn to scale) we see that $w_{\nu,m}^{-} > \beta_{m}^{-}$, noting that $\mathcal{Z}_{3}(z) > \zeta(z)$ (see (\ref{eq1.22}), (\ref{eq1.23}) and (\ref{eq1.32})), with both functions monotonically decreasing (see (\ref{eq1.35a}) and (\ref{eq1.35})). Note in the figure we have shown $w_{\nu,m}^{-}<z_{m,0}$ but this is not necessarily true or detrimental if not so.

Now, on noting that $\zeta^{-1}(\zeta_{m,0})=z_{m,0}$ (see (\ref{eq1.5})), let
\begin{equation}
\label{eqA.45}
B_{m}^{-}(v)
:=\left.\frac{\beta_{m}^{-}}{z_{m,0}} 
\right\vert_{\zeta_{m,0} = v/(v-1)}
=\frac{\zeta^{-1}(0.989\,v/(v-1))}
{\zeta^{-1}(v/(v-1))}.
\end{equation}
In \cref{fig:InverseFig} $B_{m}^{-}(v)$ is plotted for $0 \leq v < 1$ (corresponding to $-\infty < \zeta_{m.0} \leq 0$). We observe that it is bounded below by its value as $v \to 1^{-}$ ($\zeta_{m,0} \to -\infty$). We find numerically that this value is $0.98354\cdots$. We conclude that $\beta_{m}^{-} > 0.9835 \,z_{m,0}$ and since $w_{\nu,m}^{-} > \beta_{m}^{-}$ the lower bound of (\ref{eq1.60}) follows.

\begin{figure}[htbp]
 \centering
 \includegraphics[width=0.7\textwidth,
 keepaspectratio]
 {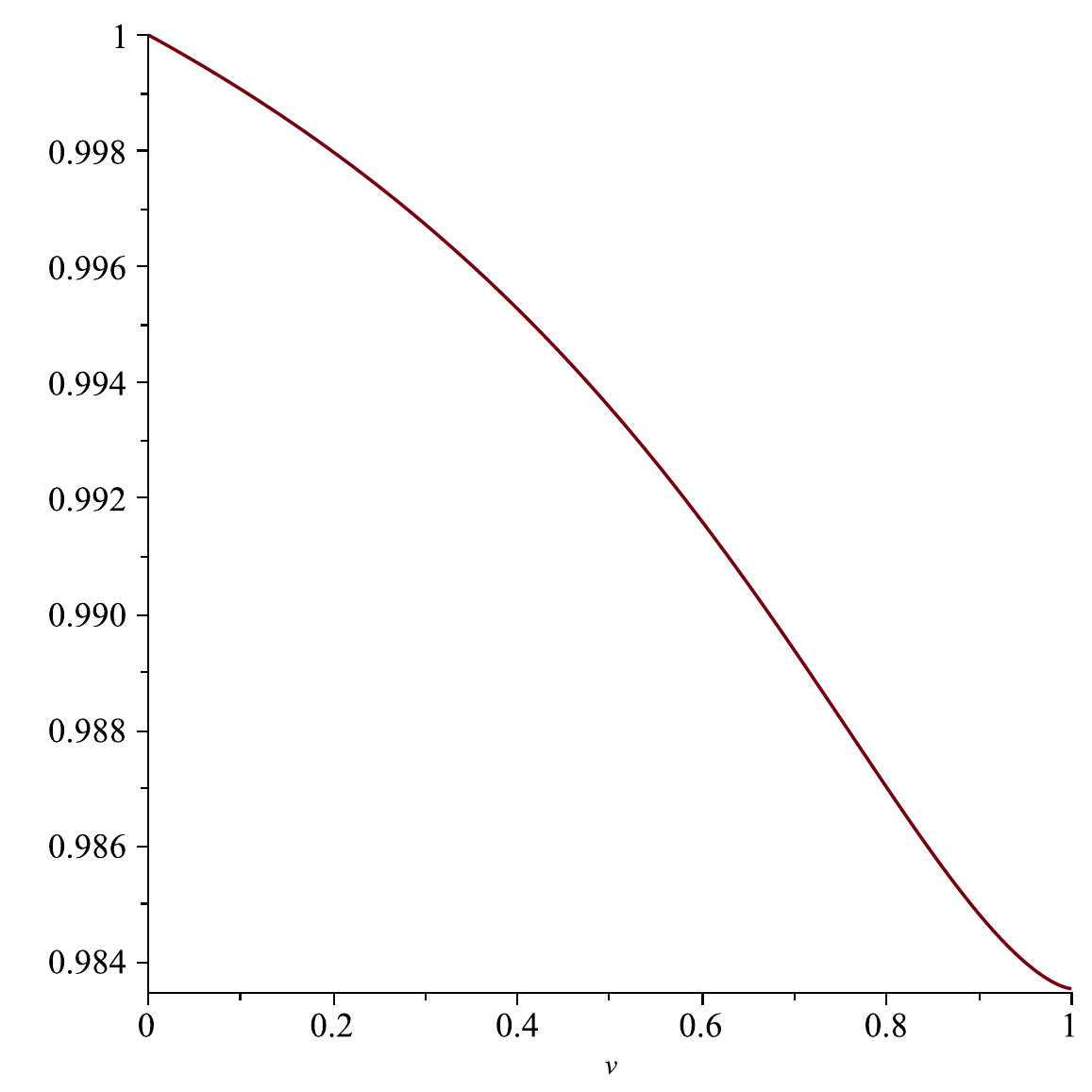}
 \caption{Graph of $B_{m}^{-}(v)$ for $0 \leq v < 1$}
 \label{fig:InverseFig}
\end{figure}

The proof of the upper bound in (\ref{eq1.60}) follows similarly. Thus define
\begin{equation}
\label{eqA.46}
\alpha_{\nu,m}^{+}
:=\mathcal{Z}_{3}(\nu,w_{\nu,m}^{+})
=\frac{\mathrm{a}_{m,0}}{\nu^{2/3}}
\left(1+\frac{0.01+0.03\delta_{m}}{4m-1}\right)
\quad  (m=1,2,3,\ldots),
\end{equation}
where we recall that $\delta_{m}=1$ for $m=1,2$ and is zero otherwise. We then have
\begin{equation}
\label{eqA.47}
\frac{\alpha_{\nu,m}^{+}}{\zeta_{m,0}}
=\frac{\mathrm{a}_{m,0}}{\mathrm{a}_{m}}
\left(1+\frac{0.01+0.03\delta_{m}}{4m-1}\right)
\quad  (m=1,2,3,\ldots).
\end{equation}
Note from (\ref{eq1.5}), (\ref{eq1.57}) and (\ref{eqA.46}), and recalling $r_{m}^{-}>0$ from \cref{thm:QuWong}, that $|\alpha_{\nu,m}^{+}|>|\zeta_{m,0}|$ for $m=1,2,3,\ldots$, $1 \leq \nu < \infty$.

An upper bound is derived as follows. Firstly, it is readily verified from explicit computation of (\ref{eqA.47}) that 
\begin{equation}
\label{eqA.47a}
\frac{\alpha_{\nu,2}^{+}}{\zeta_{2,0}}
< \frac{\alpha_{\nu,1}^{+}}{\zeta_{1,0}}
= 1.00559 \cdots < 1.0056.
\end{equation}
Next, again from \cite[Eqs. (4.4) and (4.5)]{Qu:1999:BPU} and (\ref{eqA.47}),
\begin{equation}
\label{eqA.48}
\frac{\alpha_{\nu,m}^{+}}{\zeta_{m,0}}
\leq
\left\{1-\frac{0.130}{\left(
\frac{3}{8}\pi(4m-1.051)\right)^{2}}\right\}^{-1}
\left(1+\frac{0.01}{4m-1}\right)
\quad  (m=3,4,5,\ldots),
\end{equation}
with RHS being found through elementary calculus to attain its absolute maximum at $m=3$, the value at which is $1.00169\cdots < 1.0056$. Thus we have for $m=1,2,3,\ldots$, $1 \leq \nu < \infty$
\begin{equation}
\label{eqA.49}
|\zeta_{m,0}|<|\alpha_{\nu,m}^{+}|<1.0056|\zeta_{m,0}|.
\end{equation}
Then similarly to our earlier lower bound on $w_{\nu,m}^{-}$, and again referring to \cref{fig:zeta-Z}, it follows that $w_{\nu,m}^{+} <\beta_{\nu,m}^{+}:=\mathcal{Z}_{3}^{-1}(\nu,1.0056\,\zeta_{m,0})$ (this being the inverse function with respect to the second argument).

\begin{figure}[htbp]
 \centering
 \includegraphics[width=0.7\textwidth,
 keepaspectratio]
 {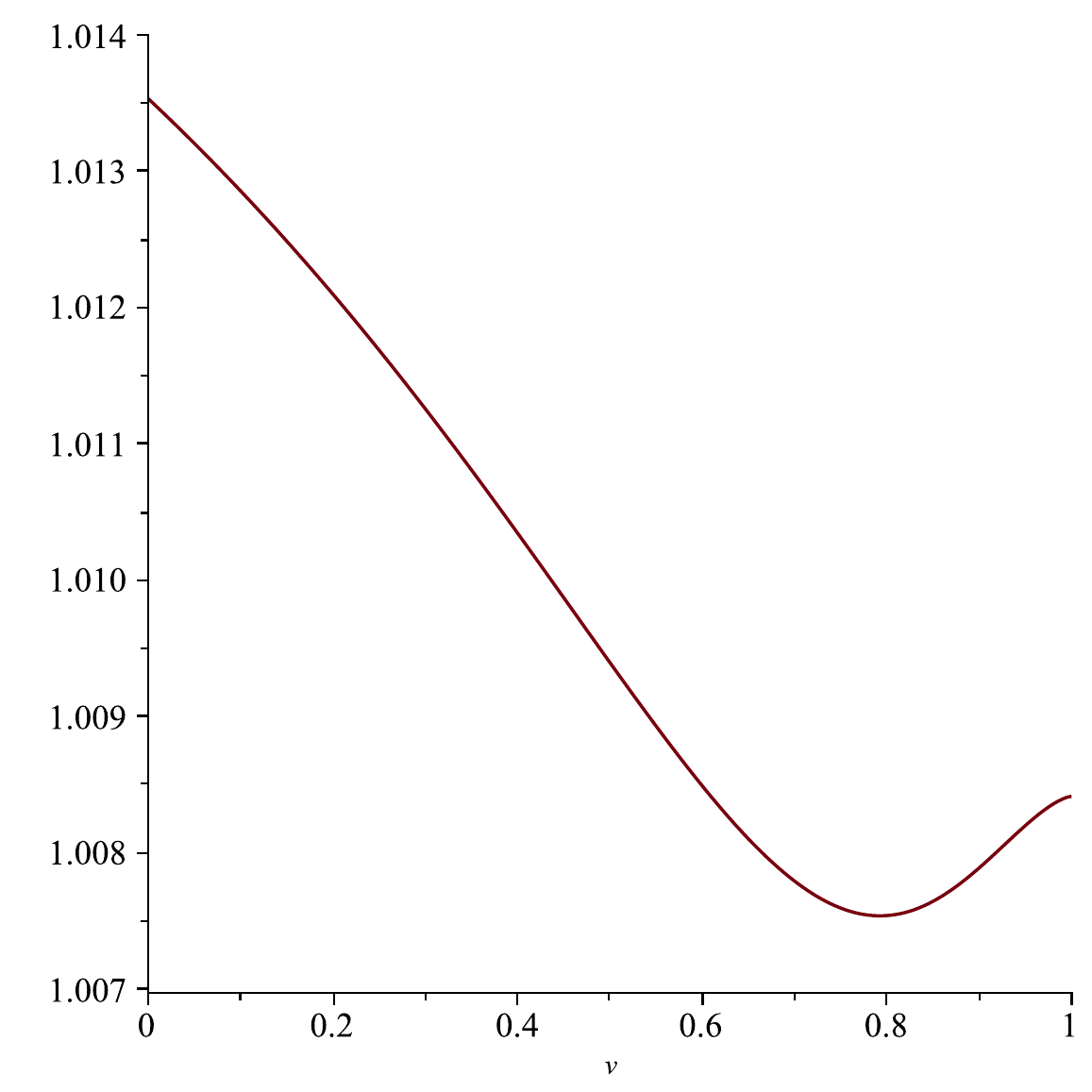}
 \caption{Graph of $B_{m}^{+}(v)$ for $0 \leq v < 1$}
 \label{fig:InverseFig2}
\end{figure}

We now use the inequality $\beta_{\nu,m}^{+} \leq \beta_{1,m}^{+}$ for $1 \leq \nu < \infty$, which follows in this case from the fact that $\mathcal{Z}_{3}(\nu,z) \leq \mathcal{Z}_{3}(1,z)$ (see \cref{lem:SumUpsilon} and (\ref{eq1.22})). With this in mind let
\begin{equation}
\label{eqA.49a}
B_{m}^{+}(v)=\left.\frac{\beta_{1,m}^{+}}{z_{m,0}} 
\right\vert_{\zeta_{m,0} = v/(v-1)}
=\frac{\mathcal{Z}_{3}^{-1}(1,1.0056\,v/(v-1))}
{\zeta^{-1}(v/(v-1))},
\end{equation}
which is plotted in \cref{fig:InverseFig2} for $0 \leq v <1$. We observe that it attains a maximum value of $1.0135313599\cdots$ at $v=0$ ($\zeta_{m,0}=0$). As a result $w_{\nu,m}^{+} <\beta_{\nu,m}^{+} \leq \beta_{1,m}^{+} < 1.01354 \, z_{m,0}$, and as such the upper bound of (\ref{eq1.60}) has been established, completing the proof of the lemma.
\end{proof}

\begin{proof}[Proof of \cref{lem:sigma}]

From (\ref{eq1.3}) and integration by parts, for $z>1$ ($\zeta<0$),
\begin{equation*}
2|\zeta|^{3/2}
=\frac{\left(z^2-1\right)^{3/2}}{z^2}
+2\int_{1}^{z}
\frac{\left(t^2-1\right)^{3/2}}{t^3}\,dt
>\frac{\left(z^2-1\right)^{3/2}}{z^2}.
\end{equation*}
Further
\begin{equation*}
2|\zeta|^{3/2}
=3 \int_{1}^{z}
\frac{\left(t^2-1\right)^{1/2}}{t}\,dt
< 3\int_{1}^{z} t\left(t^2-1\right)^{1/2}dt
= \left(z^2-1\right)^{3/2},
\end{equation*}
and
\begin{equation*}
2|\zeta|^{3/2}
=3 \int_{1}^{z}
\frac{\left(t^2-1\right)^{1/2}}{t}\,dt
< \frac{3(z-1)\left(z^2-1\right)^{1/2}}{z},
\end{equation*}
since $t^{-1}(t^2-1)^{1/2}$ is increasing for $t \geq 1$. Thus on dividing all three by $2(z^2-1)^{3/2}$ we obtain the bounds
\begin{equation}
\label{eqA.32c}
\frac{1}{2z^2} <\sigma^3
< \min\left\{\frac{3}{2 z (z+1)}, \frac{1}{2} \right\}.
\end{equation}
Thus, from (\ref{eq1.3}), (\ref{eq1.8}) and (\ref{eqA.32c}), again for $z>1$ ($\zeta<0$),
\begin{equation}
\label{eqA.31}
\frac{d \sigma}{dz}
=-\frac{2z^2\sigma^3-1}{2 z |\zeta|} <0,
\end{equation}
and from the product rule, (\ref{eqA.32c}) and (\ref{eqA.31})
\begin{equation}
\label{eqA.33}
\frac{d (z\sigma)}{dz}
=\frac{2 \sigma \left(z^2\sigma^2 + \zeta\right)
- 1}{2\zeta}
=\frac{1- 2 \sigma^3}{2|\zeta|} >0,
\end{equation}
noting from (\ref{eq1.8}) that $\zeta=-\sigma^2 (z^2-1)$. The asserted monotonicity on $\sigma(z)$ and $z\sigma(z)$ follows.
\end{proof}

\begin{proof}[Proof of \cref{lem:etaBound}]

First, using (\ref{eq1.35a}) and (\ref{eq1.34}),
\begin{equation*}
\dot{\eta} =
\frac{d \eta(\nu,z)}{d z}
\left(\frac{d\zeta}{dz}\right)^{-1}
=-z\sigma(z)\eta'(\nu,z) > 0.
\end{equation*}
Further, from \cref{lem:sigma,lem:UpsilonMonotonic} $\sigma(z)$ and $z|\eta'(1,z)|$ are decreasing for $1<z<\infty$, and thus from (\ref{eq1.9}) and (\ref{eq1.34})
\begin{equation*}
1  < 1+\dot{\eta}
\leq \sup_{1 \leq z < \infty} \left\{
1-z\sigma(z)\eta'(\nu,z) \right\} 
\leq  1-\sigma(1)\eta'(1,1) 
=\tfrac{2611707229667}{2590330640625},
\end{equation*}
and the result follows.
\end{proof}

\begin{proof}[Proof of \cref{lem:MBound}]

From (\ref{eq1.60}) we have under the hypothesis of the lemma that $z \geq w_{\nu,m}^{-} \geq \max\{0.9835z_{m,0},1\}$, and since from \cref{lem:zmjMonotonic} $\hat{z}_{3}(z)$ is positive and monotonically decreasing we deduce $\hat{z}_{3}(z) \leq \hat{z}_{3}(\max\{0.9835z_{m,0},1\})$. Now it is readily shown numerically that $p_{16}(v):=
\hat{z}_{3}(\max\{0.9835z,1\})/\hat{z}_{3}(z)$, where $z=1/(1-v)$, is monotonically increasing for $0 \leq v <1$ ($1 \leq z < \infty$); see \cref{fig:zm3rat}. Thus on referring to (\ref{eq1.19}) 
\begin{equation*}
\hat{z}_{3}(\max\{0.9835 z,1\})
/\hat{z}_{3}(z)
< \lim_{z \to \infty} \left\{\hat{z}_{3}(0.9835z)
/\hat{z}_{3}(z)\right\} = (0.9835)^{-5},
\end{equation*}
and so on setting $z=z_{m,0}$ in this, and recalling $\hat{z}_{3}(z_{m,0})=z_{m,3}$, establishes (\ref{eq1.62}).

\begin{figure}[htbp]
 \centering
 \includegraphics[width=0.7\textwidth,
 keepaspectratio]
 {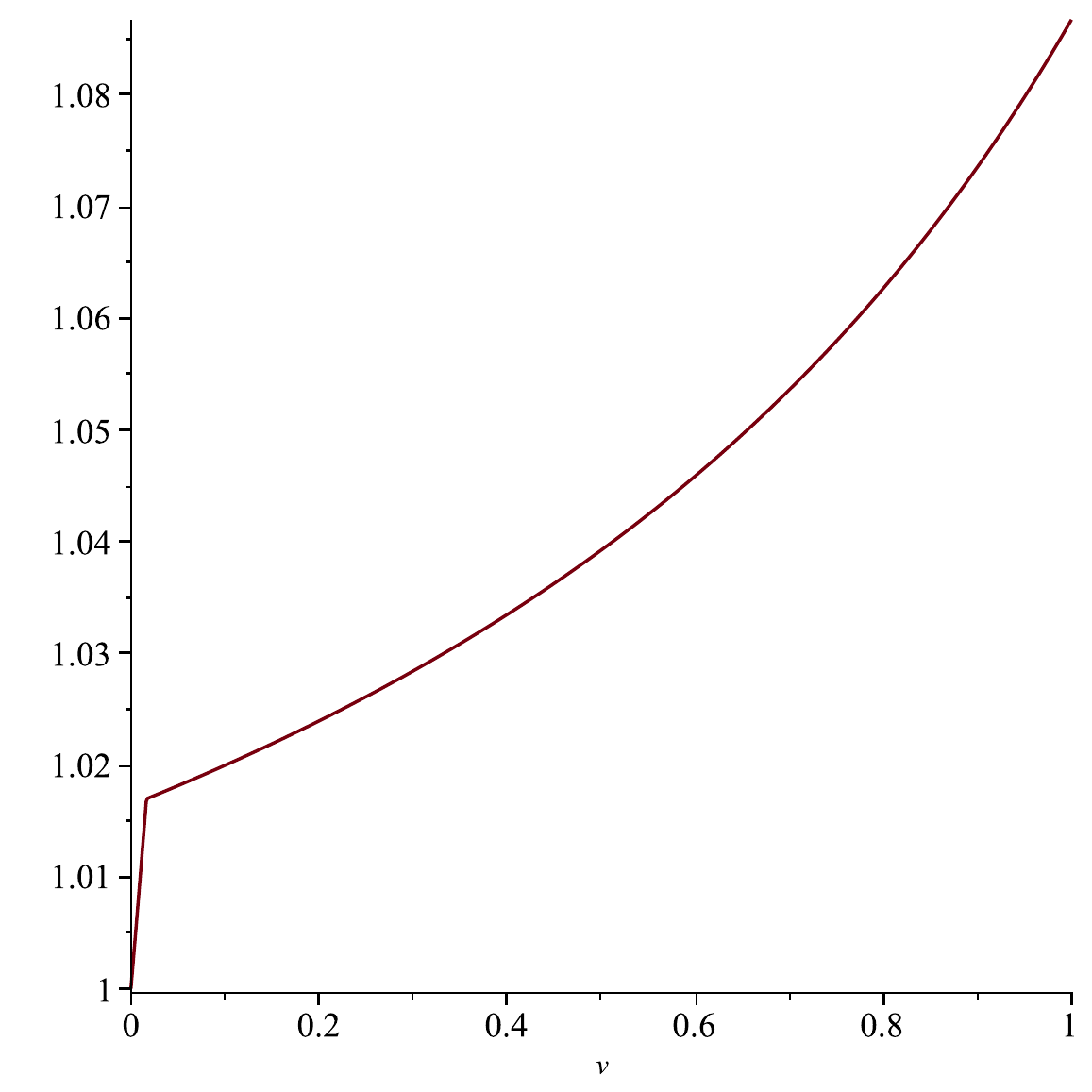}
 \caption{Graph of $p_{16}(v)$ for $0 \leq v < 1$}
 \label{fig:zm3rat}
\end{figure}

Finally, from (\ref{eq1.22}), (\ref{eq1.23}), (\ref{eq1.33}) and (\ref{eq1.57})
\begin{equation*}
\nu^{2/3}\zeta(w_{\nu,m}^{-})
=\mathrm{a}_{m}+r_{m}^{+}
-\nu^{2/3}\eta(\nu,w_{\nu,m}^{-})
<\mathrm{a}_{m}+r_{m}^{+}.
\end{equation*}
Then under the hypothesis $z \geq w_{\nu,m}^{-}$, and recalling that $w_{\nu,m}^{-}>1$, we have $\zeta=\zeta(z) \leq \zeta(w_{\nu,m}^{-})$ since $\zeta$ is decreasing for $1 \leq z < \infty$ (see (\ref{eq1.35a})). Therefore since $M(x)$ is increasing (see \cref{remark:Mincreasing}) it follows that $z \geq w_{\nu,m}^{-} \implies$
\begin{equation*}
M(\nu^{2/3}\zeta)
\leq M(\nu^{2/3}\zeta(w_{\nu,m}^{-}))
< M(\mathrm{a}_{m}+r_{m}^{+})
\leq \kappa' \pi^{-1/2}\left|\mathrm{a}_{m,0}\right|^{-1/4},
\end{equation*}
where
\begin{equation*}
\kappa'
=\sup_{m \in \mathbb{Z}^{+}}
\left\{
\sqrt{\pi}\left|\mathrm{a}_{m,0}\right|^{1/4}
M(\mathrm{a}_{m}+r_{m}^{+})\right\}
= 1.0000277286\cdots,
\end{equation*}
with the supremum being attained at $m=12$. In this computation we used (\ref{eq1.50}), (\ref{eq1.39}) and (\ref{eq1.52}). The veracity of the bound (\ref{eq1.63}) is now evident.
\end{proof}

\begin{proof}[Proof of \cref{lem:zhatbound}]

From (\ref{eq1.35}) and (\ref{eq1.43}) we observe that $\partial F_{m}(\nu,z)/ \partial z=\mathcal{Z}_{3}'(\nu,z)<0$ for $1 \leq z < \infty$ and $1 \leq \nu < \infty$, and hence $F_{m}(\nu,z)$ is a strictly decreasing function of $z$ in this interval. Furthermore, for fixed $m=1,2,3,\ldots$ it is readily verified from (\ref{eq1.3c}), (\ref{eq1.29}) - (\ref{eq1.31}) that
\begin{equation}
\label{eqA.27}
F_{m}(\nu,z) \to -\infty
\quad (z \to \infty).
\end{equation}
But from (\ref{eq1.5}) and (\ref{eq1.32}) $F_{m}(\nu,z_{m,0})=\eta(\nu,z_{m,0})>0$ ($1 \leq \nu < \infty$). Thus by the strict monotonicity of $F_{m}(\nu,z)$, and the opposite signs of both functions at the points $z=z_{m,0}$ and $z=\infty$, we have shown that $\hat{z}_{\nu,m}$ is indeed the unique simple zero in the interval $(z_{m,0},\infty)$, with the lower bound having been established.

It remains to prove the upper bound of (\ref{eq1.44}), and to this end we use the Taylor remainder theorem and (\ref{eq1.35a}) to obtain the identity
\begin{equation}
\label{eqA.29}
\zeta\left(z_{m,0}+\tfrac{1}{73}\right)
=\zeta\left(z_{m,0}\right)
+\frac{\zeta'(\tau)}{73}
=\frac{\mathrm{a}_{m}}{\nu^{2/3}}
-\frac{1}{73\tau \sigma(\tau)},
\end{equation}
for some number $\tau$ satisfying
\begin{equation}
\label{eqA.30}
z_{m,0}<\tau<z_{m,0}+\tfrac{1}{73}.
\end{equation}

Now from \cref{lem:sigma} $\{z \sigma\}^{-1}$ is strictly decreasing, and hence from (\ref{eqA.29}) and (\ref{eqA.30})
\begin{equation*}
\zeta\left(z_{m,0}+\tfrac{1}{73}\right)
<\frac{\mathrm{a}_{m}}{\nu^{2/3}}
-\frac{1}{73 (z_{m,0}+\tfrac{1}{73})
\sigma\left(z_{m,0}+\tfrac{1}{73}\right)}.
\end{equation*}
It follows from (\ref{eq1.5}), (\ref{eq1.22}), (\ref{eq1.23}), (\ref{eq1.32}), and (\ref{eq1.43}) that
\begin{equation*}
\begin{split}
F_{m}\left(\nu,z_{m,0}+\tfrac{1}{73}\right) &
< -\frac{1}{73(z_{m,0}+\tfrac{1}{73})
\sigma\left(z_{m,0}+\tfrac{1}{73}\right)}
+ \eta\left(\nu,z_{m,0}+\tfrac{1}{73}\right)
\\ &
< -\frac{1}{73(z_{m,0}+\tfrac{1}{73})
\sigma\left(z_{m,0}+\tfrac{1}{73}\right)}
+ \eta\left(1,z_{m,0}+\tfrac{1}{73}\right).
\end{split}
\end{equation*}
Consequently, since $1+\frac{1}{73} < z_{m,0}+\frac{1}{73} < \infty$ for $m=1,2,3,\ldots$ and $1 \leq \nu < \infty$,
\begin{equation}
\label{eqA.36}
F_{m}\left(\nu,z_{m,0}+\tfrac{1}{73}\right)
<-\frac{1-c_{2}}{73(z_{m,0}+\tfrac{1}{73})
\sigma\left(z_{m,0}+\tfrac{1}{73}\right)},
\end{equation}
where
\begin{equation}
\label{eqA.37}
c_{2} :=73 \sup_{\frac{74}{73} < z < \infty}
\left\{z \sigma(z) \eta(1,z)\right\}
= 74 \, \sigma\left(\tfrac{74}{73}\right)
\eta\left(1,\tfrac{74}{73}\right)
=0.99176 \cdots <1,
\end{equation}
with the supremum being attained at $z={74}/{73}$ since  $z \sigma(z) \eta(1,z)$ is positive and montonically decreasing for $1 \leq z<\infty$; see \cref{lem:UpsilonMonotonic}. Thus from (\ref{eqA.36}) and (\ref{eqA.37}) it is seen that $F_{m}(\nu,z_{m,0}+\tfrac{1}{73})<0$ for $m=1,2,3,\ldots$ and $1 \leq \nu < \infty$. Accordingly, from (\ref{eqA.27}) and the strict monotonicity of the function, we deduce that its (sole) simple zero in $(1,\infty)$, namely $z=\hat{z}_{\nu,m}$, must be smaller than $z=z_{m,0}+\tfrac{1}{73}$, and the asserted upper bound in (\ref{eq1.44}) follows.
\end{proof}

\section*{Acknowledgements}
I thank the referees for helpful comments, and acknowledge financial support from Ministerio de Ciencia e Innovación project PID2021-127252NB-I00 (MCIN/AEI/10.13039/ 501100011033/FEDER, UE).

\bibliographystyle{siamplain}
\bibliography{biblio}
\end{document}